\documentclass[12pt]{article}
\usepackage[long]{edale}
\usepackage[all]{xy}
\usepackage{mathrsfs}
\DeclareMathOperator{\Hom}{Hom}
\DeclareMathOperator{\sExt}{\underline{Ext}}
\DeclareMathOperator{\Spec}{Spec}
\DeclareMathOperator{\Spa}{Spa}
\DeclareMathOperator{\Spf}{Spf}
\DeclareMathOperator{\pr}{pr}
\DeclareMathOperator{\Fr}{Fr}
\DeclareMathOperator{\Tr}{Tr}
\DeclareMathOperator{\cl}{cl}
\DeclareMathOperator{\Aut}{Aut}
\DeclareMathOperator{\Gal}{Gal}
\DeclareMathOperator{\CH}{CH}
\DeclareMathOperator{\id}{id}
\DeclareMathOperator{\Gr}{Gr}
\DeclareMathOperator{\obj}{obj}
\DeclareMathOperator{\ch}{ch}
\DeclareMathOperator{\Sw}{Sw}
\DeclareMathOperator{\dimtot}{dimtot}
\DeclareMathOperator{\length}{length}

\newcommand{\hooklongrightarrow}{\lhook\joinrel\longrightarrow}
\newcommand{\A}{\mathbb{A}}
\newcommand{\Z}{\mathbb{Z}}
\newcommand{\F}{\mathbb{F}}
\newcommand{\Q}{\mathbb{Q}}
\newcommand{\C}{\mathbb{C}}
\newcommand{\X}{\mathsf{X}}
\newcommand{\rig}{\mathrm{rig}}
\newcommand{\adj}{\mathrm{adj}}

\newcommand{\bc}{\mathrm{b.c.}}
\newcommand{\cf}{\textit{cf.\ }}
\newcommand{\Lotimes}{\overset{\mathbb{L}}{\otimes}}
\newcommand{\Lboxtimes}{\overset{\mathbb{L}}{\boxtimes}}

\makeatletter
\renewcommand{\maketitle}%
 {\@ifundefined{@title}{\relax}%
  {%
   \vspace*{30pt}\begin{center}\Large{\bfseries\mathversion{bold}\@title}\end{center}%
   \@ifundefined{@author}{\relax}{\begin{center}\large\@author\end{center}}%
   \thispagestyle{plain}%
  }%
 }
\makeatother

\newcommand{\yleftarrow}[2][]{\settowidth{\dimen0}{\ensuremath{#1}}\settowidth{\dimen1}{\ensuremath{#2}}%
\xleftarrow[\ifthenelse{\lengthtest{\dimen0>12.5pt}}{#1}{\makebox[12.5pt]{\ensuremath{\scriptstyle #1}}}]%
{\ifthenelse{\lengthtest{\dimen1>12.5pt}}{#2}{\makebox[12.5pt]{\ensuremath{\scriptstyle #2}}}}}

\title{On $\ell$-independence for the \'etale cohomology of rigid spaces over local fields}
\author{Yoichi Mieda}
\begin{document}
\maketitle

\begin{firstfootnote}
 Graduate School of Mathematical Sciences, the University of Tokyo,
 3--8--1 Komaba, Meguro-ku, Tokyo 153--8914, Japan.

 E-mail address: \texttt{edale@ms.u-tokyo.ac.jp}

 2000 \textit{Mathematics Subject Classification}.
 Primary: 14F20;
 Secondary: 14G20, 14G22.
\end{firstfootnote}

\begin{abstract}
 We investigate the action of the Weil group on the compactly supported $\ell$-adic \'etale cohomology groups
 of rigid spaces over a local field. 
 We prove that the alternating sum of the traces of the action is an integer and is independent of $\ell$
 when either the rigid space is smooth or the characteristic of the base field is equal to $0$.
 We modify the argument of T.~Saito (\cite{TSaito}) to prove a result on $\ell$-independence
 for nearby cycle cohomology, which leads to our $\ell$-independence result for smooth rigid spaces.
 In the general case, we use the finiteness theorem of R.~Huber (\cite{Huber-JAG-!}),
 which requires the restriction on the characteristic of the base field.
\end{abstract}

\section{Introduction}
Let $K$ be a complete discrete valuation field with finite residue field $\F_q$ and
$\overline{K}$ a separable closure of $K$.
We denote by $\Fr_q$ the geometric Frobenius element (the inverse of the $q$th power map)
in $\Gal(\overline{\F}_q/\F_q)$.
The Weil group $W_K$ of $K$ is defined as the inverse image of the subgroup
$\langle \Fr_q\rangle\subset\Gal(\overline{\F}_q/\F_q)$ by the canonical map
$\Gal(\overline{K}/K)\longrightarrow \Gal(\overline{\F}_q/\F_q)$.
For $\sigma\in W_K$, let $n(\sigma)$ be the integer such that the image of $\sigma$ in $\Gal(\overline{\F}_q/\F_q)$
is $\Fr_q^{n(\sigma)}$. Put $W_K^+=\{\sigma\in W_K\mid n(\sigma)\ge 0\}$.

Let $\X$ be a separated rigid space over $K$. 
We consider the action of $W_K$ on the compactly supported $\ell$-adic cohomology group 
$H^i_c(\X\otimes_K\overline{K},\Q_\ell)$, where $\ell$ is a prime number which does not divide $q$. 
This cohomology group us defined by using the \'etale site of $\X$ (\cf \cite{Huber-book}, \cite{Huber-comparison}).
Our main theorem is the following:

\begin{thm}[][Theorem \ref{Thm:smooth-l-ind}, Theorem \ref{Thm:general-l-ind}]{Theorem}\label{l-independence}
Let $\mathsf{X}$ be a quasi-compact separated rigid space over $K$. Assume one of the following conditions:
\begin{itemize}
 \item The rigid space $\mathsf{X}$ is smooth over $K$.
 \item The characteristic of $K$ is equal to $0$.
\end{itemize}
Then for every $\sigma \in W_K^+$, the number
\[
 \sum_{i=0}^{2\dim \mathsf{X}}(-1)^i\Tr\bigl(\sigma_*;H^i_c(\mathsf{X}\otimes_K\overline{K},\Q_\ell)\bigr)
\] 
is an integer which is independent of $\ell$.
\end{thm}

Note that $H^i_c(\X\otimes_K\overline{K},\Q_\ell)$ is known to be a finite-dimensional $\Q_\ell$-vector space
when one of the above conditions is satisfied (\cite[Proposition 6.1.1, Proposition 6.2.1]{Huber-book},
\cite[Corollary 2.3]{Huber-JAG-!}, \cite[Theorem 3.1]{Huber-comparison}). 
In the previous paper, under the same assumption, the author proved that every eigenvalue of the action
of $\sigma\in W_K^+$ on $H^i_c(\X\otimes_K\overline{K},\Q_\ell)$ is a Weil number 
(\cite[Theorem 4.2, Theorem 5.5]{rigid-Weil}). 

For a scheme over $K$, the property in Theorem \ref{l-independence} was proven by Ochiai (\cite[Theorem 2.4]{Ochiai}).
However it seems difficult to prove Theorem \ref{l-independence} by the same method as in \cite{Ochiai},
since the induction on the dimension does not work well. 
In this paper, we modify the method in \cite{TSaito},
which treats the composite action of an element of $W_K$ and a correspondence.

We sketch the outline of the paper. In \S \ref{section:finite-l-ind}, we derive $\ell$-independence
of the alternating sum of the traces of the action of a correspondence from Fujiwara's trace formula
(\cite{Fujiwara-trace}). This result seems well-known, but we include its proof for completeness.
In \S \ref{section:cycle-class}, by using localized Chern characters, 
we prove a lemma which is a refined version of \cite[Lemma 2.17]{TSaito}.
This lemma is needed in \S \ref{section:analogue-of-weight-spec-seq}.
In \S \ref{section:partially-supp-coh}, we introduce partially supported cohomology and investigate
its several functorial properties. In terms of partially supported cohomology,
we can describe the action of a correspondence on the compactly supported cohomology
of a scheme which is not necessarily proper.
The required properties of nearby cycles and their cohomology are also included in this section.
In \S \ref{section:analogue-of-weight-spec-seq}, we introduce a spectral sequence converging to nearby cycle cohomology,
which is a generalization of the weight spectral sequence studied in \cite{RaZi}, \cite{TSaito}.
By the same method as in \cite[\S 2.3, \S 2.4]{TSaito}, we can prove the compatibility of the spectral
sequence with the action of a correspondence.
In \S \ref{section:nearby-l-ind}, we prove $\ell$-independence for nearby cycle cohomology
by using the result in \S \ref{section:analogue-of-weight-spec-seq} and de Jong's alteration (\cite{deJong}).
The method is almost the same as that in \cite[\S 3]{TSaito}. Several applications to algebraic geometry
(not to rigid geometry) are also included (Theorem \ref{Thm:stalk-l-ind}, Theorem \ref{Thm:open-l-ind}).
Finally in \S \ref{section:rigid-l-ind} we give a proof of our main result.
When $\X$ is smooth over $K$, we can reduce our theorem to the case where $\X$ is the Raynaud generic fiber
of the completion of a scheme over $\mathcal{O}_K$ with smooth generic fiber (though the reduction does not
seem so immediate in comparison with \cite{rigid-Weil}).
In this case we can use the result in \S \ref{section:nearby-l-ind}.
Finally, assuming that the characteristic of $K$ is $0$,
we prove our theorem for a general $\X$ by induction on $\dim \X$.
In this process, we need the finiteness theorem of R.~Huber (\cite{Huber-JAG-!}).

\bigbreak

\noindent\textbf{Acknowledgements}\ \,
The author would like to thank Takeshi Saito for valuable comments.
He would like to thank Tetsushi Ito for reading a manuscript and giving comments.
He is also grateful to his advisor Tomohide Terasoma.
He was supported by the Japan Society for the Promotion of Science Research Fellowships for
Young Scientists.

\bigbreak

\noindent\textbf{Notation}\ \,
Let $K$ be a field.
For a scheme $X$ (or a rigid space) over $K$ and an extension $L$ of $K$,
we denote the base change $X\times_{\Spec K}\Spec L$ by $X_L$.
For a scheme $X$ of finite type over $K$, we denote the group of $k$-cycles on $X$ by $Z_k(X)$ and 
the $k$th Chow group (the group of $k$-cycles modulo rational equivalences) by $\CH_k(X)$.
Let $X$ be a scheme of finite type over $K$ and $Y$ be a closed subscheme of $X$.
Put $d=\dim X$.
We denote by $\cl_Y^X\colon \CH_{d-k}(Y)\longrightarrow H^{2k}_Y(X,\Q_\ell(k))$ the cycle map defined in
\cite[cycle]{SGA4+1/2}, where $\ell$ is a prime number distinct from the characteristic of $K$.
\bigbreak

\noindent\textbf{Convention on correspondences}\ \,
Let $K$ be a field and $\ell$ a prime number distinct from the characteristic of $K$. Put $\Lambda=\Q_\ell$.
For schemes $X$ and $Y$ separated of finite type over $K$, 
a correspondence between $X$ and $Y$ is a morphism $\gamma\colon \Gamma\longrightarrow X\times Y$,
where $\Gamma$ is a scheme separated of finite type over $K$.
A morphism $f\colon X\longrightarrow X$ can be regarded as the correspondence 
$f\times \id\colon X\longrightarrow X\times X$. Note that this convention is different from that in \cite{TSaito},
while it is the same as that in \cite{SGA5} and \cite{Fujiwara-trace}.
We sometimes assume that $\gamma$ is a closed immersion.
	    
Let $\gamma\colon \Gamma\longrightarrow X\times Y$ be a correspondence such that $Y$ is smooth 
and purely $d$-dimensional. 
Put $c=\dim \Gamma$ and $\gamma_i=\pr_i\circ \gamma$. When $\gamma_1$ is proper, $\Gamma$ induces
a homomorphism between cohomology groups
\[
 \Gamma^*\colon H^q_c(X,\Lambda)\yrightarrow{\pr_1^*} H^q_c(\Gamma,\Lambda)
 \yrightarrow{\pr_{2*}} H^{q+2d-2c}_c\bigl(Y,\Lambda(d-c)\bigr).
\]
More generally, for $\alpha\in Z_k(\Gamma)$, we can define a homomorphism
\[
 \alpha^*\colon H^q_c(X,\Lambda)\longrightarrow H^{q+2d-2k}_c\bigl(Y,\Lambda(d-k)\bigr).
\]
It is easy to see that the map $\alpha^*$ depends only on the rational equivalence class of $\alpha$.
Therefore for an element $\alpha$ of the Chow group $\CH_k(\Gamma)$, we can define the map
\[
 \alpha^*\colon H^q_c(X,\Lambda)\longrightarrow H^{q+2d-2k}_c\bigl(Y,\Lambda(d-k)\bigr).
\]

\section{$\boldsymbol{\ell}$-independence for schemes over finite fields}\label{section:finite-l-ind}
\subsection{$\boldsymbol{\ell}$-independence}
\subsubsection{}
In this section we give a result on $\ell$-independence for schemes over finite fields.
Though the result seems well-known for specialists, we include its proof for completeness.

\subsubsection{Theorem}\label{Thm:finite-l-ind}
Let $X$ be a separated smooth purely $d$-dimensional scheme of finite type over $\F_q$ and 
$\gamma\colon \Gamma\longrightarrow X\times X$ a correspondence such that $\Gamma$ is purely $d$-dimensional.
We denote the characteristic of $\F_q$ by $p$.
Assume that $\gamma_1\colon \Gamma\longrightarrow X$ is proper.
Then the number
\[
 \Tr\bigl(\Gamma^*;H^*_c(X_{\overline{\F}_q},\Q_\ell)\bigr)=
 \sum_{i=0}^{2d}(-1)^i\Tr\bigl(\Gamma^*;H^i_c(X_{\overline{\F}_q},\Q_\ell)\bigr)
\]
lies in $\Z[1/p]$ and is independent of $\ell$.

\begin{prf}
 Let $\gamma^{(n)}\colon \Gamma^{(n)}\longrightarrow X\times X$ be the correspondence satisfying that
 $\gamma^{(n)}_1=\Fr_X^n\circ \gamma_1$ and $\gamma^{(n)}_2=\gamma_2$, where $\Fr_X$ is the $q$th power Frobenius
 morphism. Take a compactification $\overline{\gamma}\colon \overline{\Gamma}\longrightarrow \overline{X}\times \overline{X}$ of $\gamma\colon \Gamma\longrightarrow X\times X$ and define
 $\overline{\gamma}^{(n)}\colon \overline{\Gamma}^{(n)}\longrightarrow \overline{X}\times \overline{X}$ in the same way.
 We may assume that $D=\overline{X}\setminus X$ is a Cartier divisor of $\overline{X}$.
 Then for sufficiently large $n$, any connected component of $\overline{\Gamma}^{(n)}\cap \Delta_{\overline{X}}$
 which meets $D$ is (set-theoretically) contained in $D$ (here we identify $\Delta_{\overline{X}}$ and $\overline{X}$).
 This easily follows from Fujiwara's result on contractingness
 (\cite[Proposition 5.3.5, Proposition 5.4.1]{Fujiwara-trace}).
 See also \cite[Theorem 2.1.3, Lemma 2.2.3]{Varshavsky}.

 By this fact and Fujiwara's trace formula (\cite[Proposition 5.3.4, Proposition 5.4.1]{Fujiwara-trace}),
 there exists an integer $N$ such that for every $n\ge N$ and $\ell$ the equality
 \[
 \Tr\bigl(\Gamma^{(n)*};H^*_c(X_{\overline{\F}_q},\Q_\ell)\bigr)
 =(\Gamma^{(n)},\Delta_X)_{X\times X}
 \]
 holds. The right hand side denotes the intersection number (note that $\Gamma^{(n)}\cap \Delta_X$ is proper over
 $\F_q$ for sufficiently large $n$ by the argument above), which is an integer and is independent of $\ell$.
 Since $\Gamma^{(n)*}=\Gamma^*\circ (\Fr_X^*)^n$, the number
 $\Tr(\Gamma^*\circ (\Fr_X^*)^n;H^*_c(X_{\overline{\F}_q},\Q_\ell))$
 is an integer which is independent of $\ell$ for $n\ge N$.

 Let $\alpha_{\ell,i,1},\ldots,\alpha_{\ell,i,m_i}$ and 
 $\lambda_{\ell,i,1},\ldots,\lambda_{\ell,i,m_i}$ be eigenvalues
 of $\Gamma^*$ and $\Fr_X^*$ on $H^i_c(X_{\overline{\F}_q},\Q_\ell)$ respectively.
 By \cite[Corollaire 3.3.3, Corollaire 3.3.4]{Weil2}, $\lambda_{\ell,i,k}$ are algebraic integers and there exist
 non-negative integers $w_k$ such that the complex absolute value of any conjugate of $\lambda_{\ell,i,k}$
 is equal to $q^{w_k/2}$.
 Since $\Gamma^*$ and $\Fr_X^*$ commute with each other, the trace of $\Gamma^*\circ (\Fr_X^*)^n$ on
 $H^i_c(X_{\overline{\F}_q},\Q_\ell)$ is equal to $\sum_{k=1}^{m_i}\alpha_{\ell,i,k}\lambda_{\ell,i,k}^n$
 with $\lambda_{\ell,i,1},\ldots,\lambda_{\ell,i,m_i}$ permuted suitably. Therefore the theorem
 follows from the subsequent two lemmas.
\end{prf}

\subsubsection{Lemma}
Let $X$ be a separated $d$-dimensional scheme of finite type over $\F_q$.
Then for every eigenvalue $\lambda$ of $\Fr_X^*$ on $H^i_c(X_{\overline{\F}_q},\Q_\ell)$,
$q^d\lambda^{-1}$ is integral over $\Z$.

\begin{prf}
 We may assume that $X$ is irreducible. By de Jong's alteration \cite{deJong}, we may assume that there exist
 a proper smooth scheme $\overline{X}$ and a strict normal crossing divisor $D$ of $\overline{X}$ such that 
 $X=\overline{X}\setminus D$. Then by the Poincar\'e duality, $q^d\lambda^{-1}$ is an 
 eigenvalue of $\Fr_X^*$ on $H^i(X_{\overline{\F}_q},\Q_\ell)$.
 Let $D_1,\ldots,D_m$ be the irreducible components of $D$. Put $D_I=\bigcap_{i\in I}D_i$ for 
 $I\subset \{1,\ldots,m\}$ ($D_I=\overline{X}$ for $I=\varnothing$)
 and $D^{(k)}=\coprod_{I\subset \{1,\ldots,m\},\#I=k}D_I$.
 By the spectral sequence 
 \[
  E_1^{-k,n+k}=H^{n-k}\bigl(D^{(k)}_{\overline{\F}_q},\Q_\ell(-k)\bigr)\Longrightarrow H^n(X_{\overline{\F}_q},\Q_\ell),
 \]
 the eigenvalue $q^d\lambda^{-1}$ occurs as an eigenvalue of $\Fr_{D^{(k)}}^*$ on 
 $H^{n-k}(D^{(k)}_{\overline{\F}_q},\Q_\ell(-k))$ for some $n$, $k$.
 Since $D^{(k)}$ is proper smooth over $\F_q$,
 \cite[Corollaire 3.3.3]{Weil2} ensures that $q^d\lambda^{-1}$ is integral over $\Z$.
\end{prf}

\subsubsection{Lemma}\label{Lemma:vandermond}
Let $K$ be a field of characteristic $0$ and $\alpha_1,\ldots,\alpha_m,\lambda_1,\ldots,\lambda_m$ 
elements of $K$ such that $\lambda_k\neq 0$ for every $k$ and $\lambda_{k_1}\neq \lambda_{k_2}$ for $k_1\neq k_2$.
\begin{enumerate}
 \item Assume that there exists an integer $N$ such that $\sum_{k=1}^{m}\alpha_k\lambda_k^n=0$ for every $n\ge N$.
       Then $\alpha_1=\cdots=\alpha_m=0$.
 \item Let $q$ be a power of a prime number $p$. Assume the following conditions:
       \begin{itemize}
	\item The elements $\lambda_k$ are integral over $\Z$ and there exist
	      non-negative integers $w_k$ such that the complex absolute value of any conjugate of $\lambda_k$
	      is equal to $q^{w_k/2}$.
	\item There exists an integer $N$ such that $\sum_{k=1}^m\alpha_k\lambda_k^n\in \Z$ for every $n\ge N$.
	\item There exists an integer $d$ such that $q^d\lambda_k^{-1}$ is integral over $\Z$ for every $k$.
       \end{itemize}
       Then $\alpha_k$ is algebraic over $\Q$ and $\sum_{k=1}^m\alpha_k\lambda_k^n\in \Z[1/p]$ for every non-negative integer $n$.
\end{enumerate}

\begin{prf}
 \begin{enumerate}
  \item	This follows from
	\[
	 \det \begin{pmatrix}
	       \lambda_1^N & \cdots & \lambda_m^{N}\\
	       \vdots & & \vdots\\
	       \lambda_1^{N+m-1} & \cdots & \lambda_m^{N+m-1}
	      \end{pmatrix}
	=\lambda_1^N\cdots\lambda_m^N\prod_{i<j}(\lambda_i-\lambda_j)\neq 0.
	\]
  \item We may assume that $\alpha_k\neq 0$ for every $k$.
	Take $\sigma\in \Aut_{\Q}K$. 
	Since $\sum_{k=1}^m\sigma(\alpha_k)\sigma(\lambda_k)^n=\sum_{k=1}^m\alpha_k\lambda_k^n$ for $n\ge N$,
	by i), there exists a bijection $\tau\colon \{1,\ldots,m\}\longrightarrow \{1,\ldots,m\}$ satisfying
	$\sigma(\alpha_k)=\alpha_{\tau(k)}$ and $\sigma(\lambda_k)=\lambda_{\tau(k)}$ for every $k$.
	Especially, the set $\{\sigma(\alpha_k)\mid \sigma\in\Aut_{\Q}K\}$ is a finite set, which implies that
	$\alpha_k$ is algebraic over $\Q$. Furthermore, for every non-negative integer $n$,
	$\sigma(\sum_{k=1}^m\alpha_k\lambda_k^n)=\sum_{k=1}^m\alpha_k\lambda_k^n$. Thus
	$\sum_{k=1}^m\alpha_k\lambda_k^n\in \Q$. 

	Finally we prove that $\sum_{k=1}^m\alpha_k\lambda_k^n\in \Z[1/p]$. Let $L$ be a finite field extension of $\Q$
	generated by $\alpha_1,\ldots,\alpha_m,\lambda_1,\ldots,\lambda_m$ and $A$ the integral closure
	of $\Z$ in $L$. 
	Take a maximal ideal $\mathfrak{m}$ of $A$ which does not lie over $(p)\subset \Z$ and
	denote the valuation of $L$ associated with $\mathfrak{m}$ by $v_\mathfrak{m}$.
	By the assumption, $\lambda_k\in A_{\mathfrak{m}}^\times$.
	We should prove	$v_\mathfrak{m}(\sum_{k=1}^m\alpha_k\lambda_k^n)\ge 0$ for every $n$.
	Put $c=\max\{0,-v_{\mathfrak{m}}(\alpha_1),\ldots,-v_{\mathfrak{m}}(\alpha_m)\}$. 
	Since $A/\mathfrak{m}^c$ is a finite ring, there exists a positive integer $q$
	such that $\lambda_k^q-1\in \mathfrak{m}^c$ for every $k$.
	Take a positive integer $r$ satisfying $n+qr\ge N$. Then 
	\[
	 v_{\mathfrak{m}}\biggl(\sum_{k=1}^m\alpha_k\lambda_k^{n}\biggr)
	=v_{\mathfrak{m}}\biggl(\sum_{k=1}^m\alpha_k\lambda_k^{n+qr}
	+\sum_{k=1}^m\alpha_k(\lambda_k^{n}-\lambda_k^{n+qr})\biggr)
	\ge 0.
	\]
	This completes the proof.
 \end{enumerate}
\vspace*{-\baselineskip}
\end{prf}

\subsubsection{Remark}\label{Remark:Bloch-Esnault}
{\upshape In \cite{Bloch-Esnault}, S.~Bloch and H.~Esnault gave another proof of Theorem \ref{Thm:finite-l-ind}
by using the theory of relative motivic cohomology defined by M.~Levine.
They also prove the integrality of the alternating sum of the trace in Theorem \ref{Thm:finite-l-ind}.
They only consider the case where $X$ has a good compactification, but we can easily reduce the general case
to their case by de Jong's alteration (\cf (\ref{subsub:reduce-alteration})).
}

\section{Complements on cycle classes}\label{section:cycle-class}
\subsection{Localized Chern characters}
\subsubsection{}
Here we briefly recall localized Chern characters. Let $S$ be a noetherian regular scheme.
By an arithmetic $S$-scheme, we mean a separated scheme of finite type
over $S$. Let $\ell$ be a prime number which is invertible in $S$ and denote $\Q_\ell$ by $\Lambda$.

\subsubsection{}\label{local-chern-chow}
Let $X$ be a purely $d$-dimensional arithmetic $S$-scheme and
$i\colon Y\hooklongrightarrow X$ a closed subscheme of $X$.
Let $\mathscr{E}_\bullet$ be a bounded complex of locally free $\mathcal{O}_X$-module which is exact over
$X\setminus Y$. With such $\mathscr{E}_\bullet$, we associate $\ch^X_Y(\mathscr{E}_\bullet)\in \CH_{d-\bullet}(Y)_\Q$,
called the {\slshape localized Chern character} (\cite[18.1]{Fulton}). We denote the degree $k$-part of 
$\ch^X_Y(\mathscr{E}_\bullet)$ by $\ch^X_{k,Y}(\mathscr{E}_\bullet)\in \CH_{d-k}(Y)_\Q$. 
Note that in \cite[18.1]{Fulton}, $\ch^X_Y(\mathscr{E}_\bullet)$ is defined 
as an element of $\CH(Y\rightarrow X)_\Q$. In the notation there,
$\ch^X_Y(\mathscr{E}_\bullet)\in \CH_{d-\bullet}(Y)_\Q$ here should be denoted by 
$\ch^X_Y(\mathscr{E}_\bullet)\cap [X]$.

\subsubsection{}\label{local-chern-chow-properties}
We need the following property of $\ch^Y_X$:
\begin{quote}
 Assume that $S=\Spec K$ where $K$ is a field, $X$ is smooth over $S$ and $Y$ is irreducible. 
 Let $\mathscr{E}_\bullet\longrightarrow i_*\mathcal{O}_Y$ be a resolution of $i_*\mathcal{O}_Y$
 consisting of locally free $\mathcal{O}_X$-modules (such a resolution always exists since $X$ is regular).
 Put $d'=\dim Y$.
 Then $\ch^X_{d-d',Y}(\mathscr{E}_\bullet)=[Y]\in \CH_{d'}(Y)_{\Q}$.
 This is a corollary of the Riemann-Roch theorem (\cite[Theorem 18.3 (3), (5)]{Fulton}).
\end{quote}

\subsubsection{}\label{local-chern-etale}
Let the notation be the same as in (\ref{local-chern-chow}).
We can associate the cohomology class $\ch^X_{\ell,k,Y}(\mathscr{E}_\bullet)\in H^{2k}_Y(X,\Lambda(k))$
for each $k$, which is also called the {\slshape localized Chern character} (\cf \cite{Iversen}). 

\subsubsection{}\label{local-chern-etale-properties}
We list some properties of $\ch^X_{\ell,k,Y}$ needed later.
\begin{enumerate}
 \item The localized Chern character $\ch^X_{\ell,k,Y}(\mathscr{E}_\bullet)$ is compatible with any pull-back.
 \item Assume that $S=\Spec K$ where $K$ is a field and $X$ is smooth over $S$. 
       Then we have $\cl^X_Y(\ch^X_{k,Y}(\mathscr{E}_\bullet))=\ch^X_{\ell,k,Y}(\mathscr{E}_\bullet)$
       (\cf \cite[Example 19.2.6]{Fulton}).
\end{enumerate}

\subsection{A lemma on cycle classes}
\subsubsection{}\label{semistable-def}
Let $S=\Spec A$ be a henselian trait and $\ell$ a prime number which is invertible in $S$.
We denote the generic (resp.\ special) point of $S$ by $\eta$ (resp.\ $s$).
For an $S$-scheme $X$, we denote its generic (resp.\ special) fiber by
$X_\eta$ (resp.\ $X_s$).

An arithmetic $S$-scheme $X$ is said to be strictly semistable
if it is, Zariski locally on $X$, \'etale over $\Spec A[T_0,\ldots,T_n]/(T_0\cdots T_r-\pi)$
for a uniformizer $\pi$ of $A$ and integers $n$, $r$ with $0\le r\le n$.
Let $D_1,\ldots,D_m$ be irreducible components of $X_s$.
We put $D_I=\bigcap_{i\in I}D_i$ for $I\subset \{1,\ldots,m\}$ and $D^{(p)}=\coprod_{I\subset \{1,\ldots,m\},\#I=p+1}D_I$
for a non-negative integer $p$. We write $a_i\colon D_i\hooklongrightarrow X$ and
$a^{(p)}\colon D^{(p)}\longrightarrow X$ for the canonical morphisms.

\subsubsection{Lemma}\label{Lemma:semistable-cycle-class}
Let $X$ be a strictly semistable $S$-scheme of purely relative dimension $d$ and $Y$ a closed subscheme of $X$
with $(d-k)$-dimensional generic fiber. Assume that $Y$ is flat over $S$.
Then there exists a cohomology class $\xi_\ell\in H^{2k}_Y(X,\Q_\ell(k))$ for each prime number $\ell$
which is invertible in $S$ satisfying the following conditions:
\begin{itemize}
 \item $\xi_\ell\vert_{X_\eta}=\cl^{X_\eta}_{Y_\eta}(Y_\eta)\in H^{2k}_{Y_\eta}(X_\eta,\Q_\ell(k))$.
 \item $\xi_\ell\vert_{D^{(p)}}=\cl^{D^{(p)}}_{D^{(p)}\cap Y}(a^{(p)!}[Y])\in H^{2k}_{D^{(p)}\cap Y}(D^{(p)},\Q_\ell(k))$.
\end{itemize}
Here we are abusing notation since $D^{(p)}$ is not a subscheme of $X$.

\begin{prf}
 Take a resolution $\mathscr{E}_\bullet\longrightarrow i_*\mathcal{O}_Y$ of 
 $i_*\mathcal{O}_Y$ by locally free $\mathcal{O}_X$-modules, where $i$ denotes the canonical
 closed immersion $Y\hooklongrightarrow X$. Put $\xi_\ell=\ch^X_{\ell,k,Y}(\mathscr{E}_\bullet)$.
 Then it satisfies the first condition above by (\ref{local-chern-chow-properties}) and
 (\ref{local-chern-etale-properties}).
 
 We will prove that the second condition holds. Since the cycle map for a scheme over a field
 is compatible with the refined Gysin map, we may assume $p=0$. In other words, we should prove 
 $\xi_\ell\vert_{D_i}=\cl^{D_i}_{D_i\cap Y}(a_i^![Y])$. Since $Y$ is flat, $D_i\cap Y\hooklongrightarrow Y$ is
 a Cartier divisor. Thus $a_i^![Y]=[D_i\cap Y]$ in $\CH_{d-k}(D_i\cap Y)$. Moreover $Y$ and $D_i$ are 
 Tor-independent over $X$ and $\mathscr{E}_\bullet\vert_{D_i}$ is a resolution of $\mathcal{O}_{D_i\cap Y}$
 by locally free $\mathcal{O}_{D_i}$-modules. Therefore by (\ref{local-chern-chow-properties}) and
 (\ref{local-chern-etale-properties}), we have
 \[
  \xi_\ell\vert_{D_i}=\ch^{D_i}_{\ell,k,D_i\cap Y}(\mathscr{E}_\bullet\vert_{D_i})
 =\cl^{D_i}_{D_i\cap Y}\bigl(\ch^{D_i}_{k,D_i\cap Y}(\mathscr{E}_\bullet\vert_{D_i})\bigr)
 =\cl^{D_i}_{D_i\cap Y}(D_i\cap Y)=\cl^{D_i}_{D_i\cap Y}(a_i^![Y]).
 \]
 This completes the proof.
\end{prf}

\subsubsection{Remark}
{\upshape We can prove that the class $\xi_\ell$ constructed above coincides with the refined cycle class of $Y$
defined by using the absolute purity theorem of O.~Gabber (\cf \cite{Fujiwara-purity}).
In particular, we have the canonical element $\xi'_\ell\in H^{2k}_Y(X,\Z_\ell(k))$
whose image in $H^{2k}_Y(X,\Q_\ell(k))$ is equal to $\xi_\ell$.
}

\subsubsection{Remark}
{\upshape In the same way, we can remove the denominator $k!$ in \cite[Lemma 2.17]{TSaito}.}

\section{Partially supported cohomology and nearby cycle cohomology}\label{section:partially-supp-coh}
\subsection{Partially supported cohomology}
\subsubsection{}
Let $K$ be a separably closed field and $\ell$ a prime number which does not divide the characteristic of $K$.
Put $\Lambda=\Q_\ell$.

\subsubsection{}
Consider a triple $(X,U_1,U_2)$ of schemes over $K$ such that
\begin{quote}
 $(\star)$\quad $U_1$ is an open subscheme of $X$ and $U_2$ is an open subscheme of $U_1$.
\end{quote}
We call such a triple a {\slshape $\star$-triple}. The scheme $X$ is often assumed to be proper over $K$.
We denote the canonical open immersions $U_1\hooklongrightarrow X$ and 
$U_2\hooklongrightarrow U_1$ by $j_1$ and $j_{12}$ respectively. Put $j_2=j_1\circ j_{12}$.
A morphism $f\colon (X,U_1,U_2)\longrightarrow (Y,V_1,V_2)$ of $\star$-triples means a triple of
morphisms $f\colon X\longrightarrow Y$, $f_1\colon U_1\longrightarrow V_1$, and $f_2\colon U_2\longrightarrow V_2$
which makes the following diagram commutative:
\[
 \xymatrix{%
 U_2\ar[r]\ar[d]^{f_2}&U_1\ar[r]\ar[d]^{f_1}&X\ar[d]^{f}\\
 V_2\ar[r]&V_1\ar[r]&Y\lefteqn{.}
 }
\]

\subsubsection{Definition}
{\upshape Let $(X,U_1,U_2)$ be a $\star$-triple and $\mathcal{F}\in\obj D_c^b(U_2,\Lambda)$. 
We define the {\slshape partially supported cohomology}
$H_{!*}^q(X,U_1,U_2;\mathcal{F})$ as $H^q(X,j_{1!}Rj_{12*}\mathcal{F})$
and $H_{*!}^q(X,U_1,U_2;\mathcal{F})$ as $H^q(X,Rj_{1*}j_{12!}\mathcal{F})$.
Note that if $X$ is proper, $H_{!*}^q(X,U,U;\mathcal{F})=H^q_c(U,\mathcal{F})$.
Needless to say, $H_{*!}^q(X,U_1,U_2;\mathcal{F})=H^q(U_1,j_{12!}\mathcal{F})$ is independent of $X$.
}

\subsubsection{}\label{partial-functoriality1}
Let $f\colon (X,U_1,U_2)\longrightarrow (Y,V_1,V_2)$ be a morphism of $\star$-triples and
$k_1\colon V_1\hooklongrightarrow X$, $k_{12}\colon V_2\hooklongrightarrow V_1$ the canonical open immersions.
Put $k_2=k_1\circ k_{12}$. Consider the diagram below:
\[
 \xymatrix{%
 U_2\ar[r]^{j_{12}}\ar[d]^{f_2}&U_1\ar[r]^{j_1}\ar[d]^{f_1}&X\ar[d]^{f}\\
 V_2\ar[r]^{k_{12}}&V_1\ar[r]^{k_1}&Y\lefteqn{.}
 }
\]
Assume that one of the following conditions is fulfilled:
\begin{enumerate}
 \item The right rectangle is cartesian.
 \item The morphism $f_1$ is proper.
 \item The morphism $k_1$ is proper.
\end{enumerate}
Then we have the pull-back homomorphism
$f^*\colon H^q_{!*}(Y,V_1,V_2;\mathcal{F})\longrightarrow H^q_{!*}(X,U_1,U_2;f_2^*\mathcal{F})$
induced by the composite
\[
 k_{1!}Rk_{12*}\mathcal{F}\yrightarrow{\adj} Rf_*f^*k_{1!}Rk_{12*}\mathcal{F}\yrightarrow{\bc}
 Rf_*j_{1!}f_1^*Rk_{12*}\mathcal{F}\yrightarrow{\bc}Rf_*j_{1!}Rj_{12*}f_2^*\mathcal{F},
\]
where b.c.\ denotes the base change map.
Moreover if $f$ is proper (for example $X$ and $Y$ are proper over $K$), 
we have the push-forward homomorphism
$f_*\colon H^q_{*!}(X,U_1,U_2;Rf_2^!\mathcal{F})\longrightarrow H^q_{*!}(Y,V_1,V_2;\mathcal{F})$
induced by the composite
\[
 Rf_!Rj_{1*}j_{12!}Rf_2^!\mathcal{F}\yrightarrow{\bc}Rf_!Rj_{1*}Rf_1^!k_{12!}\mathcal{F}
 \yrightarrow{\bc}Rf_!Rf^!Rk_{1*}k_{12!}\mathcal{F}\yrightarrow{\adj}Rk_{1*}k_{12!}\mathcal{F}.
\]

\subsubsection{}\label{partial-functoriality2}
Assume that one of the following conditions is fulfilled:
\begin{enumerate}
 \item The left rectangle is cartesian.
 \item The morphism $f_2$ is proper.
 \item The morphism $k_{12}$ is proper.
\end{enumerate}
Then we have
$f^*\colon H^q_{*!}(Y,V_1,V_2;\mathcal{F})\longrightarrow H^q_{*!}(X,U_1,U_2;f_2^*\mathcal{F})$
defined similarly.
Moreover if $f$ is proper, we have
$f_*\colon H^q_{!*}(X,U_1,U_2;Rf_2^!\mathcal{F})\longrightarrow H^q_{!*}(Y,V_1,V_2;\mathcal{F})$.

\subsubsection{}
Next we define a cup product.
Let $(X,U_1,U_2)$ be a $\star$-triple such that $X$ is proper and $\mathcal{F},\mathcal{G}\in 
\obj D_{\textit{ctf}}^b(U_2,\Lambda)$.
By the lemma below, we can define a cup product
\[
 H^p_{!*}(X,U_1,U_2;\mathcal{F})\otimes H^q_{*!}(X,U_1,U_2;\mathcal{G})\yrightarrow{\cup}
 H^{p+q}_c(U_2,\mathcal{F}\Lotimes\mathcal{G})
\]
as the composite
\begin{align*}
 H^p_{!*}(X,U_1,U_2;\mathcal{F})\otimes H^q_{*!}(X,U_1,U_2;\mathcal{G})
 &=H^p(X,j_{1!}Rj_{12*}\mathcal{F})\otimes H^q(X,Rj_{1*}j_{12!}\mathcal{G})\\
 &\longrightarrow H^{p+q}(X,j_{1!}Rj_{12*}\mathcal{F}\Lotimes Rj_{1*}j_{12!}\mathcal{G})\\
 &\cong H^{p+q}\bigl(X,j_{2!}(\mathcal{F}\Lotimes \mathcal{G})\bigr)\\
 &=H^{p+q}_c(U_2,\mathcal{F}\Lotimes \mathcal{G}).
\end{align*}

\subsubsection{Lemma}
Let the notation be the same as above.
We have the isomorphism
$j_{1!}Rj_{12*}\mathcal{F}\Lotimes Rj_{1*}j_{12!}\mathcal{G}\cong j_{2!}(\mathcal{F}\Lotimes \mathcal{G})$.

\begin{prf}
 Denote the canonical closed immersion $X\setminus U_2\hooklongrightarrow X$ by $i$. 
 Since $j^*(j_{1!}Rj_{12*}\mathcal{F}\Lotimes Rj_{1*}j_{12!}\mathcal{G})=\mathcal{F}\Lotimes \mathcal{G}$,
 \[
  j_{2!}(\mathcal{F}\Lotimes \mathcal{G})\longrightarrow
 j_{1!}Rj_{12*}\mathcal{F}\Lotimes Rj_{1*}j_{12!}\mathcal{G}\longrightarrow 
 i^*(j_{1!}Rj_{12*}\mathcal{F}\Lotimes Rj_{1*}j_{12!}\mathcal{G})\yrightarrow{+1}
 \]
 is a distinguished triangle. Moreover $i^*Rj_{1*}j_{12!}\mathcal{G}=0$ implies 
 $i^*(j_{1!}Rj_{12*}\mathcal{F}\Lotimes Rj_{1*}j_{12!}\mathcal{G})=0$. Thus
 $j_{1!}Rj_{12*}\mathcal{F}\Lotimes Rj_{1*}j_{12!}\mathcal{G}\cong j_{2!}(\mathcal{F}\Lotimes \mathcal{G})$.
\end{prf}

\subsubsection{}\label{Kunneth}
Let $X$, $Y$ be proper schemes over $K$ and $U\subset X$, $V\subset Y$ open subschemes. For 
$\mathcal{F}\in\obj D_{\textit{ctf}}^b(U,\Lambda)$ and $\mathcal{G}\in\obj D_{\textit{ctf}}^b(V,\Lambda)$,
we have the following K\"unneth formula:
\[
 H^q_{!*}(X\times Y,U\times Y,U\times V;\mathcal{F}\Lboxtimes\mathcal{G})=
 H^q_{*!}(X\times Y,X\times V,U\times V;\mathcal{F}\Lboxtimes\mathcal{G})=
 \bigoplus_{i+j=q}H^i_c(U,\mathcal{F})\otimes H^j(V,\mathcal{G}).
\]

\begin{prf}
 Denote the canonical open immersions $U\hooklongrightarrow X$ and $V\hooklongrightarrow Y$ by $j$ and $k$ respectively.
 By the K\"unneth formula (\cite[Th\'eor\`eme 5.4.3]{SGA4-III}, \cite[Finitude, Th\'eor\`eme 1.9]{SGA4+1/2}), we have
 \[
  (j\times 1)_!R(1\times k)_*(\mathcal{F}\Lboxtimes\mathcal{G})=
  R(1\times k)_*(j\times 1)_!(\mathcal{F}\Lboxtimes\mathcal{G})=
 j_!\mathcal{F}\Lboxtimes Rk_*\mathcal{G}.
 \]
 This completes the proof.
\end{prf}

\subsubsection{}\label{partial-functoriality-constant1}
We write $H_{!*}^q(X,U_1,U_2)$ and $H_{*!}^q(X,U_1,U_2)$ for $H_{!*}^q(X,U_1,U_2;\Lambda)$ and 
$H_{*!}^q(X,U_1,U_2;\Lambda)$ respectively. Let $f\colon (X,U_1,U_2)\longrightarrow (Y,V_1,V_2)$ be a morphism
of $\star$-triples.
If the condition in (\ref{partial-functoriality1}) is satisfied, we have
\[
 f^*\colon H_{!*}^q(Y,V_1,V_2)\longrightarrow H_{!*}^q(X,U_1,U_2;f_2^*\Lambda)=H_{!*}^q(X,U_1,U_2).
\]
Assume further that $f$ is proper, $V_2$ is smooth and $U_2$, $V_2$ are equidimensional. Then we have
\[
 f_*\colon H_{*!}^{q+2d}(X,U_1,U_2)(d)\longrightarrow
 H_{*!}^q(X,U_1,U_2;Rf_2^!\Lambda)\longrightarrow H_{*!}^q(Y,V_1,V_2)
\]
where $d=\dim U_2-\dim V_2$. 

It is easy to see that $f^*$ and $f_*$ are dual to each other and the following projection formula holds:

\subsubsection{Proposition}\label{Prop:proj-formula1}
Assume that $X$ and $Y$ are proper over $K$.
For every $x\in H_{!*}^p(Y,V_1,V_2)$ and $y\in H_{*!}^q(X,U_1,U_2)$, the equality
$f_{2*}(f^*(x)\cup y)=x\cup f_*(y)$ holds in $H^{p+q-2d}_c(V_2,\Lambda(-d))$.

\subsubsection{}\label{partial-functoriality-constant2}
Next assume that the condition in (\ref{partial-functoriality2}) is satisfied
for a morphism $f\colon (X,U_1,U_2)\longrightarrow (Y,V_1,V_2)$ of $\star$-triples.
Then we have
\[
 f^*\colon H_{*!}^q(Y,V_1,V_2)\longrightarrow H_{*!}^q(X,U_1,U_2;f_2^*\Lambda)=H_{*!}^q(X,U_1,U_2).
\]
Assume further that $f$ is proper, $V_2$ is smooth and $U_2$, $V_2$ are equidimensional. Then we have
\[
 f_*\colon H_{!*}^{q+2d}(X,U_1,U_2)(d)\longrightarrow 
 H_{!*}^q(X,U_1,U_2;Rf_2^!\Lambda)\longrightarrow H_{!*}^q(Y,V_1,V_2)
\]
where $d=\dim U_2-\dim V_2$. 

It is easy to see that $f^*$ and $f_*$ are dual to each other and the following projection formula holds:

\subsubsection{Proposition}\label{Prop:proj-formula2}
Assume that $X$ and $Y$ are proper over $K$.
For every $x\in H_{*!}^p(Y,V_1,V_2)$ and $y\in H_{!*}^q(X,U_1,U_2)$, the equality
$f_{2*}(f^*(x)\cup y)=x\cup f_*(y)$ holds in $H^{p+q-2d}_c(V_2,\Lambda(-d))$.

\subsubsection{}\label{partial-cycle-class}
Let $(X,U_1,U_2)$ be a $\star$-triple such that $U_2$ is smooth and equidimensional.
Let $Y$ be a closed subscheme of $X$ which is purely of codimension $c$. Assume $Y\cap U_1=Y\cap U_2$ and
put $V=Y\cap U_1$. Then the diagram
\[
 \xymatrix{
 V\ar[d]^-{i_2}\ar@{=}[r]& V\ar[d]^{i_1}\ar[r]& Y\ar[d]^{i}\\
 U_2\ar[r]^-{j_{12}}& U_1\ar[r]^-{j_1}& X
 }
\]
is cartesian and we have the base change map
$Ri_2^!\Lambda=\id_!Ri_2^!\Lambda\longrightarrow Ri_1^!j_{12!}\Lambda$.
By this, we have the morphisms
$Rj_{1*}i_{1*}Ri_2^!\Lambda\longrightarrow Rj_{1*}j_{12!}\Lambda$ and 
\[
 H_V^{2c}\bigl(U_2,\Lambda(c)\bigr)\longrightarrow H_{*!}^{2c}(X,U_1,U_2)(c).
\]

\subsubsection{Lemma}\label{Lemma:partial-cycle-class}
The image of $\cl^{U_2}_V(V)\in H_V^{2c}(U_2,\Lambda(c))$ under the map above
is equal to the image of $1\in H^0(V,\Lambda)=H^0_{*!}(Y,V,V)$ under the map
$i_*\colon H^0_{*!}(Y,V,V)\longrightarrow H_{*!}^{2c}(X,U_1,U_2)(c)$.

\begin{prf}
 By the definition of $i_*$, we have the following commutative diagram:
 \[
  \xymatrix{%
 H^0(V,\Lambda)\ar@{=}[r]\ar[d]^{i_{2*}}& H^0_{*!}(Y,V,V)\ar[d]^{i_*}\\
 H^{2c}_V\bigl(U_2,\Lambda(c)\bigr)\ar[r]& H^{2c}_{*!}(X,U_1,U_2)(c)\lefteqn{,}
 }
 \]
 where the map $i_{2*}$ is induced by the canonical map $\Lambda\longrightarrow Ri_2^!\Lambda(c)[2c]$.
 By \cite[cycle, Th\'eor\`eme 2.3.8 (i)]{SGA4+1/2}, we have $i_{2*}(1)=\cl^{U_2}_V(V)$. This completes the proof.
\end{prf}

\subsubsection{}\label{partial-corr}
Let $X$, $Y$ be schemes proper over $K$ and $j\colon U\hooklongrightarrow X$, $j'\colon V\hooklongrightarrow Y$
dense open subschemes of $X$, $Y$ respectively.
Assume that $U$, $V$ are equidimensional and $V$ is smooth. Put $c=\dim U$ and $d=\dim V$.
Let $\Gamma\subset U\times V$ be a purely $d$-dimensional closed subscheme such that 
$\Gamma\hooklongrightarrow U\times V\yrightarrow{\pr_1}U$ is proper and $\overline{\Gamma}$ the closure
of $\Gamma$ in $X\times Y$.
Then $(X\times Y,U\times Y,U\times V)$ and $\overline{\Gamma}$ satisfy the condition in
(\ref{partial-cycle-class}). 
By Lemma \ref{Lemma:partial-cycle-class}, we can describe the action $\Gamma^*$ of the correspondence $\Gamma$
by means of partially supported cohomology:

\subsubsection{Proposition}\label{Prop:partial-corr}
Let the notation be the same as in (\ref{partial-corr}). 
Then $\Gamma^*\colon H^q_c(U,\Lambda)\longrightarrow H^q_c(V,\Lambda)$ coincides with the composite 
\begin{align*}
 H^q_c(U,\Lambda)&=H_{!*}^q(X,U,U)\yrightarrow{\pr_1^*}H_{!*}^q(X\times Y,U\times Y,U\times V)
 \yrightarrow{\cup \cl(\Gamma)}H_c^{q+2c}(U\times V)(c)\\
 &\yrightarrow{\pr_{2*}}H_c^q(V,\Lambda).
\end{align*}
Here $\cl(\Gamma)$ denotes the image of $\cl^{U\times V}_{\Gamma}(\Gamma)\in H^{2c}_{\Gamma}(U\times V,\Lambda(c))$ in
$H_{*!}^{2c}(X\times Y,U\times Y,U\times V)(c)$.

\begin{prf}
 This follows immediately from Proposition \ref{Prop:proj-formula1} and Lemma \ref{Lemma:partial-cycle-class}.
\end{prf}

\subsubsection{}\label{refined-pullback}
Let $f\colon X'\longrightarrow X$ be a proper morphism of equidimensional schemes over $K$ and 
$Z\subset X$ a closed subscheme of $X$ which is purely $k$-dimensional. Put $Z'=Z\times_X X'$ and $d=\dim X'-\dim X$.
Assume that $X$ is smooth over $K$. Then the map $\id\times f\colon X'\longrightarrow X'\times X$ is a regular immersion.
By applying the construction in \cite[Chapter 6]{Fulton} to the cartesian diagram
\[
 \xymatrix{Z'\ar[r]\ar[d]&X'\times Z\ar[d]\\ X'\ar[r]& X'\times X\lefteqn{,}
 }
\]
we have the element $(\id\times f)^![Z]\in \CH_{k+d}(Y')$. We denote it by $f^![Z]$.
It is well-known that this construction is compatible with cycle class, i.e., 
$f^*(\cl^X_Z(Z))=\cl^{X'}_{Z'}(f^![Z])$. 

\subsubsection{Lemma}\label{Lemma:refined-pullback-corr}
Let $X$, $Y$, $X'$ and $Y'$ be equidimensional schemes over $K$ and
$f\colon X'\longrightarrow X$, $g\colon Y'\longrightarrow Y$ be proper surjective generically finite
morphisms over $K$. Put $c=\dim X=\dim X'$ and $d=\dim Y=\dim Y'$.
Let $\Gamma\subset X\times Y$ be a purely $d$-dimensional subscheme such that the composite
$\Gamma\hooklongrightarrow X\times Y\yrightarrow{\pr_1}X$ is proper. 
Assume that $X$ and $Y$ are smooth over $K$.
Then the following diagram is commutative:
 \[
 \xymatrix{%
 H^i_c(X',\Lambda)\ar[rr]^-{((f\times g)^![\Gamma])^*}\ar[d]^-{f_*}&&
 H^i_c(Y',\Lambda)\\
 H^i_c(X,\Lambda)\ar[rr]^-{\Gamma^*}&&
 H^i_c(Y,\Lambda) \ar[u]_-{g^*}
 }
 \]

\begin{prf}
 First note that $((f\times g)^![\Gamma])^*$ makes sense since $(f\times g)^![\Gamma]$ is supported on 
 $(f\times g)^{-1}(\Gamma)$, which is proper over $X'$. 
 Take a compactification $\overline{f}\colon \overline{X}'\longrightarrow \overline{X}$ of
 $f\colon X'\longrightarrow X$ and $\overline{g}\colon \overline{Y}'\longrightarrow \overline{Y}$
 of $g\colon Y'\longrightarrow Y$. Since $f$ and $g$ are proper, we have $\overline{f}^{-1}(X)=X'$ and 
 $\overline{g}^{-1}(Y)=Y'$.
 Put $\xi=\cl(\Gamma)\in H^{2c}_{*!}(\overline{X}\times \overline{Y},X\times \overline{Y},X\times Y)(c)$. 
 Then $\cl((f\times g)^![\Gamma])=(\overline{f}\times \overline{g})^*\xi$.
 
 Consider the following diagram:
 \[
 \xymatrix{%
 H^i_c(X',\Lambda)\ar[r]^-{\pr_1^*}\ar[d]^-{f_*}
 & H^i_{!*}(\overline{X}'\times \overline{Y}',X'\times \overline{Y}',X'\times Y')\ar[rr]^-{\cup (\overline{f}\times \overline{g})^*\xi}\ar[d]^-{(\overline{f}\times \id)_*}
 && H^{i+2c}_c\bigl(X'\times Y',\Lambda(c)\bigr)\ar[r]^-{\pr_{2*}}\ar[d]^-{(f\times\id)_*}
 & H^i_c(Y',\Lambda)\\
 H^i_c(X,\Lambda)\ar[r]^-{\pr_1^*}\ar@{=}[d]
 & H^i_{!*}(\overline{X}\times \overline{Y}',X\times \overline{Y}',X\times Y')\ar[rr]^-{\cup (\id\times \overline{g})^*\xi}
 && H^{i+2c}_c\bigl(X\times Y',\Lambda(c)\bigr)\ar[r]^-{\pr_{2*}}
 & H^i_c(Y',\Lambda)\ar@{=}[u]\\
 H^i_c(X,\Lambda)\ar[r]^-{\pr_1^*}
 & H^i_{!*}(\overline{X}\times \overline{Y},X\times \overline{Y},X\times Y)\ar[rr]^-{\cup \xi}\ar[u]_-{(\id\times \overline{g})^*}
 && H^{i+2c}_c\bigl(X\times Y,\Lambda(c)\bigr)\ar[r]^-{\pr_{2*}}\ar[u]_-{(\id\times g)^*}
 & H^i_c(Y,\Lambda)\lefteqn{.}\ar[u]_-{g^*}
 }
 \]
 By Proposition \ref{Prop:partial-corr}, the composite of the upper horizontal arrows is equal to 
 $((f\times g)^![\Gamma])^*$ and that of the lower horizontal arrows is equal to $\Gamma^*$.
 The lower left rectangle, the lower middle rectangle, and the upper right rectangle in the diagram above
 are clearly commutative.
 The upper left rectangle and the lower right rectangle are commutative
 by the K\"unneth formula. The upper middle rectangle is commutative by the projection formula.
 This completes the proof.
\end{prf}

\subsection{Nearby cycle cohomology}\label{subsection:nearby-cycle-coh}
\subsubsection{}
Let $S=\Spec A$ be a strict henselian trait and denote its generic (resp.\ special) point by $\eta$ (resp.\ $s$).
Let $K$ be a quotient field of $A$ and $\overline{K}$ a separable closure of $K$. 
For an $S$-scheme $X$, we denote its special fiber, generic fiber, geometric generic fiber
by $X_s$, $X_\eta$, $X_{\overline{\eta}}$ respectively.
Denote the integral closure of $A$ in $\overline{K}$ by $\overline{A}$ and put $\overline{S}=\Spec \overline{A}$.
For an $S$-scheme $f\colon X\longrightarrow S$, we write $\overline{f}\colon \overline{X}\longrightarrow \overline{S}$
for the base change of $f$ from $S$ to $\overline{S}$. Then we have the cartesian diagrams below:
\[
 \xymatrix{%
 X_s\ar[r]^{i}\ar[d]^{f_s}& X\ar[d]^{f}& X_{\eta}\ar[l]_{j}\ar[d]^{f_\eta}\\
 s\ar[r]& S& \eta\lefteqn{,}\ar[l]
 }\qquad\qquad
 \xymatrix{%
 X_s\ar[r]^{\overline{i}}\ar[d]^{f_s}& \overline{X}\ar[d]^{\overline{f}}& X_{\overline{\eta}}\ar[l]_{\overline{j}}\ar[d]^{f_{\overline{\eta}}}\\
 s\ar[r]& \overline{S}& \overline{\eta}\lefteqn{.}\ar[l]
 }
\]

Let $\ell$ be a prime which is invertible on $S$ and denote $\Lambda=\Q_\ell$.
For $\mathcal{F}\in \obj D_c^b(X_\eta,\Lambda)$, we define 
$R\psi_X\mathcal{F}=\overline{i}^*R\overline{j}_*\varphi^*\mathcal{F}$, where 
$\varphi\colon X_{\overline{\eta}}\longrightarrow X_\eta$ is the canonical morphism.
If no confusion occurs, we omit the subscript $X$ of $R\psi_X$.

\subsubsection{}
First we recall some functorialities of nearby cycles. Let $f\colon X\longrightarrow Y$ be a morphism
between $S$-schemes. We define $f^*\colon R\psi_Y\Lambda\longrightarrow R{f_s}_*R\psi_X\Lambda$ as the composite of
\[
 R\psi_Y\Lambda\longrightarrow Rf_{s*}f_s^*R\psi_Y\Lambda\yrightarrow{\bc}
 Rf_{s*}R\psi_X f_{\overline{\eta}}^*\Lambda=Rf_{s*}R\psi_X\Lambda.
\]
Assume further that $X_\eta$, $Y_\eta$ are equidimensional and $Y_\eta$ is smooth. Put $d=\dim X_\eta-\dim Y_\eta$.
Then we define $f^*\colon Rf_{s!}R\psi_Y\Lambda(d)[2d]\longrightarrow R\psi_X\Lambda$ as the composite of
\[
 Rf_{s!}R\psi_X\Lambda(d)[2d]\yrightarrow{\bc}R\psi_Y Rf_{\overline{\eta}!}\Lambda(d)[2d]
 \longrightarrow R\psi_Y Rf_{\overline{\eta}!}Rf_{\overline{\eta}}^!\Lambda\yrightarrow{\adj}R\psi_Y\Lambda.
\]

\subsubsection{Lemma}\label{Lemma:nearby-alt}
Let $X$ and $Y$ be arithmetic $S$-schemes
and $f\colon X\longrightarrow Y$ a proper surjective generically finite $S$-morphism. 
Denote the degree of $f$ by $n$. Assume that $Y_\eta$ is smooth.
Then the composite
\[
 R\psi_Y\Lambda\yrightarrow{f^*}Rf_{s*}R\psi_X\Lambda\yrightarrow{f_*}R\psi_Y\Lambda
\]
is the multiplication by $n$.

\begin{prf}
 It is well-known that the homomorphism
 $\Lambda\longrightarrow Rf_{\overline{\eta}*}f_{\overline{\eta}}^*\Lambda\longrightarrow \Lambda$
 between constant sheaves on $Y_{\overline{\eta}}$ is the multiplication by $n$.
 Thus we have only to prove that the given map is equal to the composite
 \[
 R\psi_Y\Lambda\yrightarrow{R\psi_Y(\adj)} R\psi_Y Rf_{\overline{\eta}*}f_{\overline{\eta}}^*\Lambda
 \yrightarrow{R\psi_Y(\adj)} R\psi_Y\Lambda.
 \]
 Now we recall some basic facts on the base change map. For a commutative diagram
 \[
  \xymatrix{%
 Y'\ar[r]^-{g'}\ar[d]^{f'}& X'\ar[d]^{f}\\ Y\ar[r]^-{g}& X\lefteqn{,} 
 }
 \]
 the following hold:
 \begin{itemize}
  \item[(a)] For $\mathcal{F}\in\obj D_c^b(Y,\Lambda)$, the composite
	\[
	 Rg_*\mathcal{F}\yrightarrow{\adj}Rf_*f^*Rg_*\mathcal{F}\yrightarrow{Rf_*(\bc)}
	 Rf_*Rg'_*{f'}^*\mathcal{F}=Rg_*Rf'_*f'^*\mathcal{F}
	\]
	is equal to $Rg_*(\adj)$.
  \item[(b)] For $\mathcal{F}\in\obj D_c^b(X,\Lambda)$, the composite
	\[
	 g^*\mathcal{F}\yrightarrow{g^*(\adj)}g^*Rf_*f^*\mathcal{F}\yrightarrow{\bc}
	 Rf'_*g'^*f^*\mathcal{F}=Rf'_*f'^*g^*\mathcal{F}
	\]
	is equal to $\adj$.
 \end{itemize}
 (a) is nothing but the definition of the base change map. (b) is also easy.

\bigbreak
Consider the following cartesian diagram:
\[
 \xymatrix{%
 X_s\ar[r]^-{\overline{i}'}\ar[d]^{f_s}& \overline{X}\ar[d]^{\overline{f}}&
 X_{\overline{\eta}}\ar[l]_-{\overline{j}'}\ar[d]^{f_{\overline{\eta}}}\\
 Y_s\ar[r]^-{\overline{i}}& \overline{Y}& Y_{\overline{\eta}}\ar[l]_-{\overline{j}}\lefteqn{.}
 }
\]
By (b), the composite
\[
 \overline{i}^*R{\overline{j}}_*\Lambda\yrightarrow{\adj}Rf_{s*}f_s^*\overline{i}^*R\overline{j}_*\Lambda
 =Rf_{s*}\overline{i}'^*\overline{f}^*R\overline{j}_*\Lambda\yleftarrow[\cong]{\bc}
 \overline{i}^*R\overline{f}_*\overline{f}^*R\overline{j}_*\Lambda
\]
 is equal to $\overline{i}^*(\adj)$. Together with (a), we may infer that the composite
 \begin{align*}
 \overline{i}^*R{\overline{j}}_*\Lambda&\yrightarrow{\adj}Rf_{s*}f_s^*\overline{i}^*R\overline{j}_*\Lambda
 =Rf_{s*}\overline{i}'^*\overline{f}^*R\overline{j}_*\Lambda\yleftarrow[\cong]{\bc}
 \overline{i}^*R\overline{f}_*\overline{f}^*R\overline{j}_*\Lambda
 \yrightarrow{\bc}\overline{i}^*R\overline{f}_*R\overline{j}'_*f_{\overline{\eta}}^*\Lambda\\
 &=\overline{i}^*R\overline{j}_*Rf_{\overline{\eta}*}f_{\overline{\eta}}^*\Lambda
\end{align*}
 is equal to $\overline{i}^*R\overline{j}_*(\adj)$. Our claim easily follows from this.
\end{prf}

\subsubsection{}\label{nearby-corr}
Let $X$ and $Y$ be arithmetic $S$-schemes with equidimensional smooth generic fibers.
Put $c=\dim X_\eta$ and $d=\dim Y_\eta$. Let $\Gamma\subset X\times_S Y$ be a closed subscheme 
with purely $d$-dimensional generic fiber such that
$\Gamma\hooklongrightarrow X\times_S Y\yrightarrow{\pr_1}X$ is proper.
Denote the closed immersion $\Gamma\hooklongrightarrow X\times Y$ by $\gamma$ and put $\gamma_i=\pr_i\circ \gamma$.
Then we have the maps
\[
 H^q_c(X_s,R\psi_X\Lambda)\yrightarrow{\gamma_1^*}H^q_c(\Gamma_s,R\psi_{\Gamma}\Lambda),\quad
 H^q_c(\Gamma_s,R\psi_{\Gamma}\Lambda)\yrightarrow{\gamma_{2*}}H^q_c(Y_s,R\psi_Y\Lambda).
\]
We define $\Gamma^*\colon H^q_c(X_s,R\psi_X\Lambda)\longrightarrow H^q_c(Y_s,R\psi_Y\Lambda)$ as the composite
of the maps above.

\subsubsection{}
As in the previous subsection, we can consider a $\star$-triple $(X,U_1,U_2)$ over $S$.
We will always assume that $X$, $U_1$, $U_2$ are arithmetic $S$-schemes.
For such a $\star$-triple, we consider the partially supported nearby cycle cohomology
$H^q_{!*}(X_s,U_{1s},U_{2s};R\psi_{U_2}\Lambda)$ and $H^q_{*!}(X_s,U_{1s},U_{2s};R\psi_{U_2}\Lambda)$.
We have the same functorialities as in (\ref{partial-functoriality-constant1}), 
(\ref{partial-functoriality-constant2}),
the projection formulas, and the K\"unneth formula (\cf \cite[Th\'eor\`eme 4.2, Th\'eor\`eme 4.7]{Illusie-autour}).

\subsubsection{}\label{partial-nearby-cycle}
Let $(X,U_1,U_2)$ be a $\star$-triple over $S$ such that $U_{2\eta}$ is smooth and equidimensional.
Let $Y$ be a closed subscheme of $X$ such that $Y_\eta\subset X_\eta$ is purely of codimension $c$.
Assume $Y\cap U_1=Y\cap U_2$ and put $V=Y\cap U_1$.
Then as in (\ref{partial-cycle-class}), we have the canonical map
\[
 H_{V_\eta}^{2c}\bigl(U_{2\eta},\Lambda(c)\bigr)\longrightarrow
 H_{V_{\overline{\eta}}}^{2c}\bigl(U_{2\overline{\eta}},\Lambda(c)\bigr)
 \longrightarrow H_{V_s}^{2c}\bigl(U_{2s},R\psi_{U_2}\Lambda(c)\bigr)
 \longrightarrow H_{*!}^{2c}\bigl(X_s,U_{1s},U_{2s};R\psi_{U_2}\Lambda(c)\bigr).
\]

\subsubsection{Lemma}\label{Lemma:partial-nearby-cycle}
The image of $\cl^{U_{2\eta}}_{V_\eta}(V_\eta)$ by the map above coincides with the image of 
$1\in H^0(V_s,R\psi_V\Lambda)=H^0_{*!}(X_s,V_s,V_s;R\psi_V\Lambda)$ under the push-forward map
\[
 H^0_{*!}(X_s,V_s,V_s;R\psi_V\Lambda)\longrightarrow  H_{*!}^{2c}\bigl(X_s,U_{1s},U_{2s};R\psi_{U_2}\Lambda(c)\bigr).
\]
We denote it by $\cl(V_\eta)$.

\begin{prf}
 Consider the diagram below:
 \[
  \xymatrix{%
 H^0(V_{\eta},\Lambda)\ar[r]\ar[d]&
 H^0(V_{\overline{\eta}},\Lambda)\ar[r]\ar[d]& H^0(V_s,R\psi_V\Lambda)\ar@{=}[r]\ar[d]& 
 H^0_{*!}(Y_s,V_s,V_s;R\psi_V\Lambda)\ar[d]\\
 H^{2c}_{V_{\eta}}\bigl(U_{2\eta},\Lambda(c)\bigr)\ar[r]&
 H^{2c}_{V_{\overline{\eta}}}\bigl(U_{2\overline{\eta}},\Lambda(c)\bigr)\ar[r]&
 H^{2c}_{V_s}\bigl(U_{2s},R\psi_{U_2}\Lambda(c)\bigr)\ar[r]& 
 H^{2c}_{*!}\bigl(X_s,U_{1s},U_{2s};R\psi_V\Lambda(c)\bigr)\lefteqn{.}
 }
 \]
 The two left rectangles are clearly commutative. As in the proof of Lemma \ref{Lemma:partial-cycle-class}, 
 we can see that the right one is commutative. Since the image of $1\in H^0(V_{\eta},\Lambda)$ under the map 
 $H^0(V_{\eta},\Lambda)\longrightarrow H^{2c}_{V_{\eta}}\bigl(U_{2\eta},\Lambda(c)\bigr)$
 is $\cl^{U_{2\eta}}_{V_\eta}(V_{\eta})$, the lemma follows.
\end{prf}

\subsubsection{Proposition}\label{Prop:partial-nearby-corr}
Let $X$, $Y$ be proper arithmetic $S$-schemes and $U\subset X$, $V\subset Y$ open subschemes
which are arithmetic $S$-schemes. Assume that $U_\eta$ and $V_\eta$ are equidimensional and smooth,
and put $c=\dim U_\eta$, $d=\dim V_\eta$.
Let $\Gamma\subset U\times_S V$ be a closed subscheme with purely $d$-dimensional generic fiber 
such that the composite $\Gamma\hooklongrightarrow U\times_S V\yrightarrow{\pr_1}U$ is proper.
Then $\Gamma^*\colon H^q_c(U_s,R\psi_U\Lambda)\longrightarrow H^q_c(V_s,R\psi_V\Lambda)$ coincides with the composite 
\begin{align*}
 H^q_c(U_s,R\psi_U\Lambda)&=H_{!*}^q(X_s,U_s,U_s;R\psi_U\Lambda)
 \yrightarrow{\pr_1^*}H_{!*}^q(X_s\times Y_s,U_s\times Y_s,U_s\times V_s;R\psi_{U\times V}\Lambda)\\
 &\yrightarrow{\cup \cl(\Gamma_\eta)}H_c^{q+2c}\bigl(U_s\times V_s,R\psi_{U\times V}\Lambda(c)\bigr)
 \yrightarrow{\pr_{2*}}H_c^q(V_s,R\psi_V\Lambda).
\end{align*}
In particular, $\Gamma^*$ depends only on $\Gamma_\eta$ 
(as long as $\Gamma\hooklongrightarrow U\times_S V\yrightarrow{\pr_1}U$ is proper).

\begin{prf}
 This follows immediately from Lemma \ref{Lemma:partial-nearby-cycle} and the projection formula.
\end{prf}

\subsubsection{Lemma}\label{Lemma:nearby-refined-pullback-corr}
Let $X$, $Y$, $X'$ and $Y'$ be arithmetic $S$-schemes with smooth equidimensional generic fibers
and $f\colon X'\longrightarrow X$, $g\colon Y'\longrightarrow Y$ be proper surjective generically finite
$S$-morphisms. Put $c=\dim X_\eta=\dim X'_\eta$ and $d=\dim Y_\eta=\dim Y'_\eta$.
Let $\Gamma\subset X\times_S Y$ be a closed subscheme with purely $d$-dimensional generic fiber
such that the composite $\Gamma\hooklongrightarrow X\times_S Y\yrightarrow{\pr_1}X$ is proper. 
As in (\ref{refined-pullback}), 
we have $(f_\eta\times g_\eta)^![\Gamma_\eta]\in \CH_d((f_\eta\times g_\eta)^{-1}(\Gamma_\eta))$.
Take $\Gamma'\in Z_d((f\times g)^{-1}(\Gamma))$ whose image in $\CH_d((f_\eta\times g_\eta)^{-1}(\Gamma_\eta))$
is equal to $(f_\eta\times g_\eta)^![\Gamma_\eta]$. Such $\Gamma'$ is not unique, but $\Gamma'^*$ is independent
of the choice of $\Gamma'$ by the previous proposition. 
Then the following diagram is commutative:
 \[
 \xymatrix{%
 H^i_c(X'_s,R\psi_{X'}\Lambda)\ar[r]^-{\Gamma'^*}\ar[d]^-{f_*}&
 H^i_c(Y'_s,R\psi_{Y'}\Lambda)\\
 H^i_c(X_s,R\psi_X\Lambda)\ar[r]^-{\Gamma^*}&
 H^i_c(Y_s,R\psi_Y\Lambda) \ar[u]_-{g^*}\lefteqn{.}
 }
 \]

\begin{prf}
 As in the proof of Lemma \ref{Lemma:refined-pullback-corr}, we derive the commutativity from the projection formula,
the K\"unneth formula, and Proposition \ref{Prop:partial-nearby-corr}.
\end{prf}

\subsubsection{Lemma}\label{Lemma:special-nearby-cycle-class}
Let $(X,U_1,U_2)$ be a $\star$-triple over $S$ where $U_{2\eta}$ is smooth
and $i\colon Y\hooklongrightarrow X$ a closed subscheme such that $Y_\eta\hooklongrightarrow X_\eta$ is
purely of codimension $c$.
Assume that $Y\cap U_1=Y\cap U_2$ and put $V=Y\cap U_1$. 
Let $\xi\in H^{2c}_V(U_2,\Lambda(c))$ be an element satisfying $\xi\vert_{U_{2\eta}}=\cl^{U_{2\eta}}_{V_\eta}(V_\eta)$.
Then $\cl(V_\eta)\in H^{2c}_{*!}(X_s,U_{1s},U_{2s};R\psi_{U_2}\Lambda(c))$ coincides with the image
of $\xi$ under the map
\[
 H^{2c}_V\bigl(U_2,\Lambda(c)\bigr)\longrightarrow H^{2c}_{V_s}\bigl(U_{2s},\Lambda(c)\bigr)
 \longrightarrow H^{2c}_{*!}(X_s,U_{1s},U_{2s})(c)
 \longrightarrow H^{2c}_{*!}\bigl(X_s,U_{1s},U_{2s};R\psi_{U_2}\Lambda(c)\bigr).
\]

\begin{prf}
 This follows from the commutative diagram below:
 \[
  \xymatrix{%
 H^{2c}_V\bigl(U_2,\Lambda(c)\bigr)\ar[r]\ar[d]&
 H^{2c}_{V_s}\bigl(U_{2s},\Lambda(c)\bigr)\ar[r]\ar[dd]&
 H^{2c}_{*!}\bigl(X_s,U_{1s},U_{2s};\Lambda(c)\bigr)\ar[dd]\\
 H^{2c}_{V_\eta}\bigl(U_{2\eta},\Lambda(c)\bigr)\ar[d] &&\\
 H^{2c}_{V_{\overline{\eta}}}\bigl(U_{2\overline{\eta}},\Lambda(c)\bigr)\ar[r]&
 H^{2c}_{V_s}\bigl(U_{2s},R\psi_{U_2}\Lambda(c)\bigr)\ar[r]&
 H^{2c}_{*!}\bigl(X_s,U_{1s},U_{2s};R\psi_{U_2}\Lambda(c)\bigr)\lefteqn{.}
 }
 \]
\end{prf}

\subsubsection{Corollary}\label{Cor:special-nearby-corr}
Let the notation be the same as in Proposition \ref{Prop:partial-nearby-corr}.
Let $\xi\in H^{2c}_{\Gamma}(U\times_SV,\Lambda(c))$ be an element satisfying 
$\xi\vert_{U_\eta\times V_\eta}=\cl(\Gamma_\eta)$.
We denote by $\xi'$ the image of $\xi$ under the map
$H^{2c}_{\Gamma}(U\times V,\Lambda(c))\longrightarrow
 H^{2c}_{\Gamma_s}(U_s\times V_s,\Lambda(c))\longrightarrow
 H^{2c}_{*!}(X_s\times Y_s,U_s\times Y_s,U_s\times V_s)(c)$.
Then $\Gamma^*\colon H^q_c(U_s,R\psi_U\Lambda)\longrightarrow H^q_c(V_s,R\psi_V\Lambda)$ coincides with the composite 
\begin{align*}
 H^q_c(U_s,R\psi_U\Lambda)&=H_{!*}^q(X_s,U_s,U_s;R\psi_U\Lambda)
 \yrightarrow{\pr_1^*}H_{!*}^q(X_s\times Y_s,U_s\times Y_s,U_s\times V_s;R\psi_{U\times V}\Lambda)\\
 &\yrightarrow{\cup \xi'}H_c^{q+2c}\bigl(U_s\times V_s,R\psi_{U\times V}\Lambda(c)\bigr)
 \yrightarrow{\pr_{2*}}H_c^q(V_s,R\psi_V\Lambda).
\end{align*}

\begin{prf}
 Clear from Proposition \ref{Prop:partial-nearby-corr} and Lemma \ref{Lemma:special-nearby-cycle-class}.
\end{prf}

\section{An analogue of the weight spectral sequence and its functorialities}\label{section:analogue-of-weight-spec-seq}
\subsection{An analogue of the weight spectral sequence}
\subsubsection{}\label{semistable-notation}
Let $S=\Spec A$ be a strict henselian trait as in \ref{subsection:nearby-cycle-coh}.
Let $(X,U_1,U_2)$ a $\star$-triple over $S$ such that $U_2$ is strictly semistable (\cf (\ref{semistable-def})) 
over $S$. We say that such a $\star$-triple itself is strictly semistable.
We denote the irreducible components of $U_{2s}$ by $D''_1,\ldots,D''_m$.
We write $D_i$ (resp.\ $D_i'$) for the closure of $D''_i$ in $X$ (resp.\ $U_1$).
They form $\star$-triples $(D_i,D_i',D_i'')$. We have the following maps between $\star$-triples:
\[
 \xymatrix{%
 D_i''\ar[r]^{k'_i}\ar[d]^{a_i''}& D_i'\ar[r]^{k_i}\ar[d]^{a_i'}& D_i\ar[d]^{a_i}\\
 U_{2s}\ar[r]^{j_{12}}& U_{1s}\ar[r]^{j_1}& X_s\lefteqn{.}
 }
\]
For a subset $I\subset \{1,\ldots,m\}$, put $D_I=\bigcap_{i\in I}D_i$, $D'_I=\bigcap_{i\in I}D'_i$, and
$D''_I=\bigcap_{i\in I}D''_i$. For every $I$, $D''_I$ is smooth over $s$.
We write $a_I$, $a'_I$, $a''_I$, $k_I$, and $k'_I$ for the maps induced by
$a_i$, $a'_i$, $a''_i$, $k_i$ and $k'_i$ respectively.
For an integer $p$, put $D^{(p)}=\coprod_{I\subset \{1,\ldots,m\},\# I=p+1}D_I$,
$D'^{(p)}=\coprod_{I\subset \{1,\ldots,m\},\# I=p+1}D'_I$, and
$D''^{(p)}=\coprod_{I\subset \{1,\ldots,m\},\# I=p+1}D''_I$.
If $X$ is purely of relative dimension $n$ over $S$, they are purely of relative dimension $n-p$ over $s$.
We write $a^{(p)}$, $a'^{(p)}$, $a''^{(p)}$, $k^{(p)}$, and $k'^{(p)}$ for the maps induced by
$a_I$, $a'_I$, $a''_I$, $k_I$ and $k'_I$ respectively.
We have the following maps between $\star$-triples:
\[
 \xymatrix{%
 D_I''\ar[r]^{k'_I}\ar[d]^{a_I''}& D_I'\ar[r]^{k_I}\ar[d]^{a_I'}& D_I\ar[d]^{a_I}\\
 U_{2s}\ar[r]^{j_{12}}& U_{1s}\ar[r]^{j_1}& X_s\lefteqn{,}
 }\qquad\qquad
 \xymatrix{%
 D''^{(p)}\ar[r]^{k'^{(p)}}\ar[d]^{a''^{(p)}}& D'^{(p)}\ar[r]^{k^{(p)}}\ar[d]^{a'^{(p)}}& D^{(p)}\ar[d]^{a^{(p)}}\\
 U_{2s}\ar[r]^{j_{12}}& U_{1s}\ar[r]^{j_1}& X_s\lefteqn{.}
 }
\]

\subsubsection{}
By \cite[\S 2.1]{TSaito}, we have the monodromy filtration $M_\bullet$ on $R\psi_{U_2}\Lambda$.
This is a filtration in the category of perverse sheaves on $X_s$.
The filtration $M_\bullet$ defines a quasi-filtered object
 $(R\psi_{U_2}\Lambda,(M_sR\psi_{U_2}\Lambda/M_rR\psi_{U_2}\Lambda)_{s\ge r})$
of the category $D_c^b({U_2}_s,\Lambda)$ (\cite[5.2.17]{MSaito}). Since the functors $Rj_{12*}$ and $j_{1!}$ preserve
distinguished triangles, we have a quasi-filtered object 
\[
 \bigl(j_{1!}Rj_{12*}R\psi_{U_2}\Lambda,(j_{1!}Rj_{12*}(M_sR\psi_{U_2}\Lambda/M_rR\psi_{U_2}\Lambda))_{s\ge r}\bigr)
\]
of the category $D_c^b(X_s,\Lambda)$.

\subsubsection{Theorem}
Let the notation be the same as above.
The above quasi-filtered object induces the spectral sequence
\[
 E_1^{p,q}=\bigoplus_{i\ge\max(0,-p)}H^{q-2i}_{!*}(D^{(p+2i)},D'^{(p+2i)},D''^{(p+2i)})(-i)\Longrightarrow 
 H^{p+q}_{!*}(X_s,U_{1s},U_{2s};R\psi_{U_2}\Lambda).
\]

\begin{prf}
 By \cite[Lemme 5.2.18]{MSaito}, we have the spectral sequence
 \[
 E_1^{p,q}=H^{p+q}(X_s,j_{1!}Rj_{12*}\Gr_{-p}^M R\psi_{U_2}\Lambda)\Longrightarrow 
 H^{p+q}_{!*}(X_s,U_{1s},U_{2s};R\psi_{U_2}\Lambda).
 \]
 On the other hand, by \cite[Proposition 2.7]{TSaito} we have the canonical isomorphism
 \[
  j_{1!}Rj_{12*}\Bigl(\bigoplus_{p-q=r}a_*''^{(p+q)}\Lambda(-p)[-p-q]\Bigr)
 \yrightarrow{\cong}j_{1!}Rj_{12*}\Gr_r^M R\psi_{U_2}\Lambda.
 \]
 Since 
 $j_{1!}Rj_{12*}a''_{i*}\Lambda=j_{1!}a'_{i*}Rk'_{i*}\Lambda=j_{1!}a'_{i!}Rk'_{i*}\Lambda
 =a_{i!}k_{i!}Rk'_{i*}\Lambda=a_{i*}k_{i!}Rk'_{i*}\Lambda$, 
 we have the canonical isomorphism
 \begin{align*}
  H^{p+q}(X_s,j_{1!}Rj_{12*}\Gr_{-p}^M R\psi_{U_2}\Lambda)&\cong
  H^{p+q}\biggl(X_s,j_{1!}Rj_{12*}\Bigl(\bigoplus_{i\ge \max(0,-p)}a_*''^{(p+2i)}\Lambda(-i)[-p-2i]\Bigr)\biggr)\\
  &=\bigoplus_{i\ge \max(0,-p)}H^{q-2i}\bigl(X_s,a_*^{(p+2i)}k_!^{(p+2i)}Rk_*'^{(p+2i)}\Lambda(-i)\bigr)\\
  &=\bigoplus_{i\ge \max(0,-p)}H^{q-2i}_{!*}(D^{(p+2i)},D'^{(p+2i)},D''^{(p+2i)})(-i).
 \end{align*}
 This completes the proof.
\end{prf}

\subsubsection{}
If $X=U_1=U_2$ and $X$ is proper over $S$, the spectral sequence above coincides with the weight spectral sequence
in \cite{RaZi} up to sign (see \cite[p.~613]{TSaito}).

\subsection{Functoriality: pull-back}
\subsubsection{}
Let $(X,U_1,U_2)$ and $(Y,V_1,V_2)$ be strictly semistable $\star$-triples and
$f\colon (X,U_1,U_2)\longrightarrow (Y,V_1,V_2)$ a morphism of $\star$-triples.
Assume that the diagram
\[
 \xymatrix{
 U_1\ar[r]\ar[d]& X\ar[d]\\V_1\ar[r]& Y
 }
\]
is cartesian. Then we have the pull-back map 
$f^*\colon H^q_{!*}(Y,V_1,V_2;R\psi_{V_2}\Lambda)\longrightarrow H^q_{!*}(X,U_1,U_2;R\psi_{U_2}\Lambda)$.

\subsubsection{}
Let $E''_1,\ldots,E''_{m'}$ be the irreducible components of $V_{2s}$ and $E_i$ (resp.\ $E'_i$) a
closure of $E_i$ in $Y$ (resp.\ $V_1$).
As in (\ref{semistable-notation}), we have the diagram
\[
 \xymatrix{%
 E_i''\ar[r]^{l'_i}\ar[d]^{b_i''}& E_i'\ar[r]^{l_i}\ar[d]^{b_i'}& E_i\ar[d]^{b_i}\\
 V_{2s}\ar[r]^{j'_{12}}& V_{1s}\ar[r]^{j'_1}& Y_s\lefteqn{.}
 }
\]
We also define $E_I$, $E'_I$, $E''_I$, $E^{(p)}$, $E'^{(p)}$ and $E''^{(p)}$ as in (\ref{semistable-notation}).

Since $U_2$ and $V_2$ are strictly semistable, we have $f_2^*(\sum_{i=1}^{m'}E''_i)=\sum_{i=1}^mD''_i$ as
Cartier divisors on $U_2$. Therefore there exists a unique map 
$\varphi\colon \{1,\ldots,m\}\longrightarrow \{1,\ldots,m'\}$
satisfying $f_2(D''_i)\subset E''_{\varphi(i)}$ for every $i\in \{1,\ldots,m\}$.
Renumbering $D''_i$'s if necessary, we may assume that $\varphi$ is increasing.
Then we have $f(D_i)\subset E_{\varphi(i)}$ and $f_1(D'_i)\subset E'_{\varphi(i)}$ for every $i$. 
Moreover the right rectangle of the following commutative diagram is cartesian:
\[
 \xymatrix{%
 D''_i\ar[d]^{f_2}\ar[r]^{k'}& D'_i\ar[d]^{f_1}\ar[r]^{k}& D_i\ar[d]^{f}\\
 E''_{\varphi(i)}\ar[r]^{l'}& E'_{\varphi(i)}\ar[r]^{l}& E_{\varphi(i)}\lefteqn{.}
 }
\]

\subsubsection{}
For a non-negative integer $p$, we put $\mathcal{I}_{f,p}=\{I\subset \{1,\ldots,m\}\mid \# I=\# \varphi(I)=p+1\}$
and $D''^{(p)}_f=\coprod_{I\in \mathcal{I}_{f,p}}D''_I$. We define $D_f^{(p)}$ and $D'^{(p)}_f$ similarly.
For $I\in \mathcal{I}_{f,p}$, we have a morphism of $\star$-triples 
$f_{\varphi(I)I}\colon (D_I,D'_I,D''_I)\longrightarrow (E_{\varphi(I)},E'_{\varphi(I)},E''_{\varphi(I)})$,
which is a restriction of $f$. Put
$f^{(p)}=\coprod f_{\varphi(I)I}\colon (D^{(p)}_f,D'^{(p)}_f,D''^{(p)}_f)\longrightarrow (E^{(p)},E'^{(p)},E''^{(p)})$.
Since the right rectangle of the commutative diagram
\[
 \xymatrix{%
 D''^{(p)}_f\ar[d]^{f_2^{(p)}}\ar[r]^{k'^{(p)}}& D'^{(p)}_f\ar[d]^{f_1^{(p)}}\ar[r]^{k^{(p)}}&
 D^{(p)}_f\ar[d]^{f^{(p)}}\\
 E''^{(p)}\ar[r]^{l'^{(p)}}& E'^{(p)}\ar[r]^{l^{(p)}}& E^{(p)}
 }
\]
is cartesian, we have the pull-back map 
\[
 f^{(p)*}=\sum_{I\in\mathcal{I}_{f,p}}f^*_{\varphi(I)I}\colon 
 H^q_{!*}(E^{(p)},E'^{(p)},E''^{(p)})\longrightarrow H^q_{!*}(D^{(p)}_f,D'^{(p)}_f,D''^{(p)}_f)
 \hooklongrightarrow H^q_{!*}(D^{(p)},D'^{(p)},D''^{(p)}).
\]

\subsubsection{Proposition}\label{Prop:functoriality-pull-back}
We have a map of spectral sequences
\[
 \xymatrix{%
 E'^{p,q}_1=\bigoplus_{i\ge\max(0,-p)}H^{q-2i}_{!*}(E^{(p+2i)},E'^{(p+2i)},E''^{(p+2i)})(-i)
 \ar@{=>}[r]\ar[d]^-{\oplus f^{(p+2i)*}}&
 H^{p+q}_{!*}(Y_s,V_{1s},V_{2s};R\psi_{V_2}\Lambda)\ar[d]^-{f^*}\\
 E^{p,q}_1=\bigoplus_{i\ge\max(0,-p)}H^{q-2i}_{!*}(D^{(p+2i)},D'^{(p+2i)},D''^{(p+2i)})(-i)
 \ar@{=>}[r]&
 H^{p+q}_{!*}(X_s,U_{1s},U_{2s};R\psi_{U_2}\Lambda)\lefteqn{.}
 }
\]

\begin{prf}
 We have a morphism of quasi-filtered objects
 \begin{align*}
  &\bigl(j'_{1!}Rj'_{12*}R\psi_Y\Lambda,(j'_{1!}Rj'_{12*}(M_sR\psi_{V_2}\Lambda/M_rR\psi_{V_2}\Lambda))_{s\ge r}\bigr)\\
  &\quad\longrightarrow 
  \bigl(Rf_{s*}f_s^*j'_{1!}Rj'_{12*}R\psi_{V_2}\Lambda,(Rf_{s*}f_s^*j'_{1!}Rj'_{12*}(M_sR\psi_{V_2}\Lambda/M_rR\psi_{V_2}\Lambda))_{s\ge r}\bigr)\\
 &\quad\longrightarrow 
 \bigl(Rf_{s*}j_{1!}Rj_{12*}f_{2s}^*R\psi_{V_2}\Lambda,(Rf_{s*}j_{1!}Rj_{12*}f_{2s}^*(M_sR\psi_{V_2}\Lambda/M_rR\psi_{V_2}\Lambda))_{s\ge r}\bigr).
 \end{align*}
 By \cite[Proposition 2.11 (1)]{TSaito}, we have a morphism of quasi-filtered objects
 \[
  \bigl(f_{2s}^*R\psi_{V_2}\Lambda,f_{2s}^*(M_sR\psi_{V_2}\Lambda/M_rR\psi_{V_2}\Lambda)_{s\ge r}\bigr)
 \longrightarrow \bigl(R\psi_{U_2}\Lambda,(M_sR\psi_{U_2}\Lambda/M_rR\psi_{U_2}\Lambda)_{s\ge r}\bigr),
 \]
 which induces
 \begin{align*}
 &\bigl(Rf_{s*}j_{1!}Rj_{12*}f_{2s}^*R\psi_{V_2}\Lambda,Rf_{s*}j_{1!}Rj_{12*}f_{2s}^*(M_sR\psi_{V_2}\Lambda/M_rR\psi_{V_2}\Lambda)_{s\ge r}\bigr)\\
 &\quad\longrightarrow \bigl(Rf_{s*}j_{1!}Rj_{12*}R\psi_{U_2}\Lambda,Rf_{s*}j_{1!}Rj_{12*}(M_sR\psi_{U_2}\Lambda/M_rR\psi_{U_2}\Lambda)_{s\ge r}\bigr).
 \end{align*}
Therefore we have a morphism of quasi-filtered objects
 \begin{align*}
  &\bigl(j'_{1!}Rj'_{12*}R\psi_{V_2}\Lambda,(j'_{1!}Rj'_{12*}(M_sR\psi_{V_2}\Lambda/M_rR\psi_{V_2}\Lambda))_{s\ge r}\bigr)\\
  &\quad\longrightarrow \bigl(Rf_{s*}j_{1!}Rj_{12*}R\psi_{U_2}\Lambda,Rf_{s*}j_{1!}Rj_{12*}(M_sR\psi_{U_2}\Lambda/M_rR\psi_{U_2}\Lambda)_{s\ge r}\bigr).
 \end{align*}
 The associated morphism of spectral sequences is as follows:
 \[
  \xymatrix{%
 E'^{p,q}_1=H^{p+q}(Y_s,j'_{1!}Rj'_{12*}\Gr_{-p}^MR\psi_{V_2}\Lambda)\ar@{=>}[r]\ar[d]&
 H^{p+q}_{!*}(Y_s,V_{1s},V_{2s};R\psi_{V_2}\Lambda)\ar[d]^{f^*}\\
 E^{p,q}_1=H^{p+q}(X_s,j_{1!}Rj_{12*}\Gr_{-p}^MR\psi_{U_2}\Lambda)\ar@{=>}[r]&
 H^{p+q}_{!*}(X_s,U_{1s},U_{2s};R\psi_{U_2}\Lambda)\lefteqn{.}
 }
 \]
 On the other hand, by \cite[Proposition 2.11 (2)]{TSaito}, we have the following commutative diagram for every $r$:
 \[
  \xymatrix{%
 \bigoplus_{p-q=r}j'_{1!}Rj'_{12*}b_*''^{(p+q)}\Lambda(-p)[-p-q]\ar[r]^-{\cong}\ar[d]&
 j'_{1!}Rj'_{12*}\Gr_r^MR\psi_{V_2}\Lambda\ar[d]\\
 \bigoplus_{p-q=r}Rf_{s*}j_{1!}Rj_{12*}f_{2s}^*b_*''^{(p+q)}\Lambda(-p)[-p-q]\ar[r]^-{\cong}\ar[d]&
 Rf_{s*}j_{1!}Rj_{12*}f_{2s}^*\Gr_r^MR\psi_{V_2}\Lambda\ar[d]\\
 \bigoplus_{p-q=r}Rf_{s*}j_{1!}Rj_{12*}a_*''^{(p+q)}\Lambda(-p)[-p-q]\ar[r]^-{\cong}&
 Rf_{s*}j_{1!}Rj_{12*}\Gr_r^MR\psi_{U_2}\Lambda\lefteqn{,}
 }
 \]
 where the horizontal arrows are the canonical isomorphisms in \cite[Proposition 2.7]{TSaito}.
 We know that $j'_{1!}Rj'_{12*}b_*''^{(p+q)}\Lambda=b_*^{(p+q)}l_!Rl'_*\Lambda$ and
 $Rf_{s*}j_{1!}Rj_{12*}a_*''^{(p+q)}\Lambda=Rf_{s*}a_*^{(p+q)}k_!Rk'_*\Lambda$.
 Thus we have the map 
 $b^{(p+q)}_*l_!Rl'_*\Lambda\longrightarrow Rf_{s*}a^{(p+q)}_*k_!Rk'_*\Lambda$ induced by
 the composite of 
 \[
  j'_{1!}Rj'_{12*}b_*''^{(p+q)}\Lambda\longrightarrow Rf_{s*}j_{1!}Rj_{12*}f_{2s}^*b_*''^{(p+q)}\Lambda
 \longrightarrow Rf_{s*}j_{1!}Rj_{12*}a_*''^{(p+q)}\Lambda,
 \]
 appearing in the above diagram. We can easily see that (by taking $R\Gamma(Y_s,*)$) this map induces
 \[
 f^{(p+q)*}\colon H^k_{*!}(E^{(p+q)},E'^{(p+q)},E''^{p+q})\longrightarrow H^k_{*!}(D^{(p+q)},D'^{(p+q)},D''^{(p+q)}).
 \]
 The proposition immediately follows from this.
\end{prf}

\subsection{Functoriality: cup product}
\subsubsection{Proposition}\label{Prop:functoriality-cup-product}
Let $(X,U_1,U_2)$ be a strictly semistable $\star$-triple over $S$ and $\xi\in H^m_{*!}(X_s,U_{1s},U_{2s})(l)$.
Then the cup product with $\xi$ induces a map of spectral sequences
\[
 \xymatrix{%
 E_1^{p,q}=\bigoplus_{i\ge\max(0,-p)}H^{q-2i}_{!*}(D^{(p+2i)},D'^{(p+2i)},D''^{(p+2i)})(-i)
 \ar[d]^{\cup \xi\vert_{D^{(p+2i)}}}\ar@{=>}[r]&
 H^{p+q}_{!*}(X_s,U_{1s},U_{2s};R\psi_{U_2}\Lambda)\ar[d]^{\cup \xi}\\
  E_1^{p,q+m}=\bigoplus_{i\ge\max(0,-p)}H^{q-2i+m}_c\bigl(D''^{(p+2i)},\Lambda(-i+l)\bigr)\ar@{=>}[r]&
 H^{p+q+m}_c\bigl(U_{2s},R\psi_{U_2}\Lambda(l)\bigr)\lefteqn{.}
 }
\]

\begin{prf}
 We have a map of quasi-filtered objects
 \begin{align*}
  &\bigl(j_{1!}Rj_{12*}R\psi_{U_2}\Lambda,(j_{1!}Rj_{12*}(M_sR\psi_{U_2}\Lambda/M_rR\psi_{U_2}\Lambda))_{s\ge r}\bigr)\\
  &\quad \yrightarrow{\cup\xi} \bigl(j_{2!}R\psi_{U_2}\Lambda(l)[m],(j_{2!}(M_sR\psi_{U_2}\Lambda/M_rR\psi_{U_2}\Lambda)(l)[m])_{s\ge r}\bigr)
 \end{align*}
 and the map of spectral sequences induced by it:
 \[
  \xymatrix{%
 E_1^{p,q}=H^{p+q}(X_s,j_{1!}Rj_{12*}\Gr_{-p}^MR\psi_{U_2}\Lambda)\ar@{=>}[r]\ar[d]^{\cup\xi}&
 H^{p+q}_{!*}(X_s,U_{1s},U_{2s};R\psi_{U_2}\Lambda)\ar[d]^{\cup \xi}\\
  E_1^{p,q+m}=H^{p+q+m}\bigl(X_s,j_{2!}\Gr_{-p}^MR\psi_{U_2}\Lambda(l)\bigr)\ar@{=>}[r]&
 H^{p+q+m}_c\bigl(U_{2s},R\psi_{U_2}\Lambda(l)\bigr)\lefteqn{.}
 }
 \]
 By the subsequent lemma, the following diagram is commutative for every $r$:
 \[
  \xymatrix{%
 \bigoplus_{p-q=r}j_{1!}Rj_{12*}a_*''^{(p+q)}\Lambda(-p)[-p-q]\ar[r]^-{\cong}\ar[d]^{\cup\xi}& j_{1!}Rj_{12*}\Gr_r^MR\psi_{U_2}\Lambda\ar[d]^{\cup \xi}\\
 \bigoplus_{p-q=r}j_{2!}a_*''^{(p+q)}\Lambda(-p+l)[-p-q+m]\ar[r]^-{\cong}& j_{2!}\Gr_r^MR\psi_{U_2}\Lambda(l)[m]\lefteqn{.}
 }
 \]
 On the other hand, the diagram below is obviously commutative:
 \[
  \xymatrix{%
 j_{1!}Rj_{12*}a_*''^{(p+q)}\Lambda\ar@{=}[r]\ar[d]^{\cup \xi}& a_*^{(p+q)}k_!Rk'_*\Lambda\ar[d]^{\cup \xi\vert_{D^{(p+q)}}}\\
 j_{1!}Rj_{12*}a_*''^{(p+q)}\Lambda(l)[m]\ar@{=}[r]& a_*^{(p+q)}(k\circ k')_!\Lambda(l)[m]\lefteqn{.}
 }
 \]
 This completes the proof.
\end{prf}

\subsubsection{Lemma}
Let $X$ be a scheme over a field and $j\colon U\hooklongrightarrow X$ be an open subscheme.
Let $\mathcal{F}$ and $\mathcal{G}$ be objects of $D_c^b(U,\Lambda)$.
Then for every morphism $\varphi\colon \mathcal{F}\longrightarrow \mathcal{G}$ and
every cohomology class $\xi\in H^m_c(X,\Lambda(l))=\Hom(\Lambda,j_!\Lambda(l)[m])$, the following diagram is commutative:
\[
 \xymatrix{%
 Rj_*\mathcal{F}\ar[r]^{\varphi}\ar[d]^{\cup \xi}&Rj_*\mathcal{G}\ar[d]^{\cup \xi}\\
 j_!\mathcal{F}(l)[m]\ar[r]^{\varphi}&j_!\mathcal{G}(l)[m]\lefteqn{.}
 }
\]

\begin{prf}
 This follows from the diagram below, whose rectangles are easily seen to be commutative:
 \[
  \xymatrix{%
 Rj_*\mathcal{F}\ar[r]^-{\varphi}\ar[d]^-{\id\otimes\xi}& Rj_*\mathcal{G}\ar[d]^-{\id\otimes\xi}\\
 Rj_*\mathcal{F}\Lotimes j_!\Lambda(l)[m]\ar[r]& Rj_*\mathcal{G}\Lotimes j_!\Lambda(l)[m]\\
 j_!j^*\bigl(Rj_*\mathcal{F}\Lotimes j_!\Lambda(l)[m]\bigr)\ar[r]\ar[u]_{\cong}\ar[d]^{\cong}&
 j_!j^*\bigl(Rj_*\mathcal{G}\Lotimes j_!\Lambda(l)[m]\bigr)\ar[u]_{\cong}\ar[d]^{\cong}\\
 j_!(\mathcal{F}\Lotimes \Lambda(l)[m]\bigr)\ar[r]\ar[d]^{\cong}&
 j_!(\mathcal{G}\Lotimes \Lambda(l)[m]\bigr)\ar[d]^{\cong}\\
 j_!\mathcal{F}(l)[m]\ar[r]^-{\varphi}& j_!\mathcal{G}(l)[m]\lefteqn{.}
 }
 \]
\end{prf}

\subsection{Functoriality: push-forward}
\subsubsection{}
Let $X$ and $Y$ be strictly semistable $S$-schemes and $f\colon X\longrightarrow Y$ a morphism between them.
We assume that $X$ (resp.\ $Y$) be purely of relative dimension $n$
(resp.\ $n'$). Put $d=n-n'$. We denote the irreducible component of $X$ (resp.\ $Y$) by $D_1,\ldots,D_m$ 
(resp.\ $E_1,\ldots,E_{m'}$) and define $D^{(p)}$ (resp.\ $E^{(p)}$) as in (\ref{semistable-notation}).

\subsubsection{Proposition}\label{Prop:functoriality-push-forward}
We have a map of spectral sequences
\[
 \xymatrix{%
 E_1^{p,q+2d}=\bigoplus_{i\ge \max(0,-p)}H_c^{q+2d-2i}\bigl(D^{(p+2i)},\Lambda(-i+d)\bigr)
 \ar@{=>}[r]\ar[d]^{\oplus f_*^{(p+2i)}}& H^{p+q+2d}_c\bigl(X_s,R\psi_X\Lambda(d)\bigr)\ar[d]^{f_*}\\
 E'^{p,q}_1=\bigoplus_{i\ge \max(0,-p)}H_c^{q-2i}\bigl(E^{(p+2i)},\Lambda(-i)\bigr)
  \ar@{=>}[r]& H^{p+q}_c(Y_s,R\psi_Y\Lambda)\lefteqn{,}
 }
\] 
where $f^{(p)}_*$ is defined as in \cite[\S 2.3]{TSaito}.

\begin{prf}
 This follows immediately from \cite[Proposition 2.13]{TSaito}.
\end{prf}

\subsection{Functoriality: action of correspondence}
\subsubsection{}
Let $X$ and $Y$ be strictly semistable $S$-schemes and $X\hooklongrightarrow \overline{X}$, 
$Y\hooklongrightarrow \overline{Y}$ compactifications over $S$.
Assume that $X$ (resp.\ $Y$) is purely of relative dimension $n$ (resp.\ $n'$).
Let $D_1,\ldots,D_m$ (resp.\ $E_1,\ldots,E_{m'}$) be the irreducible components of $X_s$ (resp.\ $Y_s$).
Denote $\overline{D}_i$ (resp.\ $\overline{E}_i$) the closure of $D_i$ in $\overline{X}$ 
(resp.\ of $E_i$ in $\overline{Y}$). Write $\mathcal{I}_i$ (resp.\ $\mathcal{I}_i'$) for the defining ideal of
$\overline{D}_i$ (resp.\ $\overline{E}_i$). Let $\pi\colon \overline{Z}\longrightarrow \overline{X}\times_S\overline{Y}$
be the blow-up of $\overline{X}\times_S\overline{Y}$ by the ideal
$\prod_{(i,i')\in \Delta}(\prod_{j=1}^i\pr_1^*\mathcal{I}_j+\prod_{j'=1}^{i'}\pr_2^*\mathcal{I}'_{j'})$,
where $\Delta$ denotes the set $\{1,\ldots,m\}\times \{1,\ldots,m'\}$. Put $Z=\pi^{-1}(X\times Y)$.
Then by \cite[Lemma 1.9]{TSaito}, $Z$ is strictly semistable over $S$ and the irreducible components of $Z_s$
are indexed by $\Delta$ as $\{C_{i,i'}\}_{(i,i')\in \Delta}$ so that $\pi(C_{i,i'})=D_i\times E_{i'}$.
For $I''\subset \{1,\ldots,m\}\times \{1,\ldots,m'\}$, put $C_{I''}=\bigcap_{(i,i')\in I''}C_{i,i'}$.
We know that $C^{(p)}=\coprod_{\# I''=p+1}C_{I''}$ (\cite[Lemma 1.9]{TSaito}),
where $I''$ runs over all totally ordered subsets of $\Delta$ (the order of $\Delta$ is the product order).

For $I\subset \{1,\ldots,m\}$ and $I'\subset \{1,\ldots,m'\}$ satisfying $\# I=\# I'=p+1$,
denote by $I\wedge I'\subset \Delta$ the graph of the increasing bijection $I\longrightarrow I'$.
Put $C_1^{(p)}=\coprod_{I\subset \{1,\ldots,m\}, I'\subset \{1,\ldots,m'\}, \#I=\#I' =p+1}C_{I\wedge I'}$.
Let $\pi_{I\wedge I'}\colon C_{I\wedge I'}\longrightarrow D_I\times E_{I'}$ be the restriction of $\pi$
and $\pi^{(p)}\colon C_1^{(p)}\longrightarrow D^{(p)}\times E^{(p)}$ the morphism induced by $\pi_{I\wedge I'}$.

\subsubsection{}\label{functoriality-corr-notation}
Let $\Gamma\subset X\times_S Y$ be a closed subscheme with purely $n'$-dimensional generic fiber
such that the composite $\Gamma\hooklongrightarrow X\times_S Y\yrightarrow{\pr_1}X$ is proper.
Denote by $\Gamma'$ the closure of $\Gamma_\eta\subset X_\eta\times Y_\eta=Z_\eta$ in $Z$ and
put $\Gamma'^{(p)}\in \CH_{n'-p}(C^{(p)}\cap \Gamma')$ the refined pull-back of $\Gamma'$ to $C^{(p)}$.
By Lemma \ref{Lemma:semistable-cycle-class}, there exists a cohomology class 
$\xi\in H^{2n}_{\Gamma'}(X\times_S Y,\Lambda(n))$ satisfying the following conditions:
\begin{itemize}
 \item $\xi\vert_{X_\eta\times Y_\eta}=\cl_{X_\eta\times Y_\eta}(\Gamma_\eta)$.
 \item $\xi\vert_{C^{(p)}}=\cl^{C^{(p)}}_{C^{(p)}\cap \Gamma'}(\Gamma'^{(p)})$.
\end{itemize}
Since $\Gamma'\subset \pi^{-1}(\Gamma)$, the composite of
$C_1^{(p)}\cap \Gamma'\hooklongrightarrow C_1^{(p)}\yrightarrow{\pi^{(p)}}
D^{(p)}\times E^{(p)}\yrightarrow{\pr_1} D^{(p)}$
is proper. Thus $\Gamma'^{(p)}$ induces the action on cohomology
$(\Gamma'^{(p)})^*\colon H_c^q(D^{(p)},\Lambda)\longrightarrow H_c^q(E^{(p)},\Lambda)$
(we write $\Gamma'^{(p)}$ again for the restriction of $\Gamma'^{(p)}$ to $C^{(p)}_1\cap \Gamma'$).

On the other hand we have 
$\Gamma''^{(p)}=\pi^{(p)}_*(\Gamma'^{(p)})\in \CH_{n'-p}((D^{(p)}\times E^{(p)})\cap \Gamma)$.
As the composite 
$(D^{(p)}\times E^{(p)})\cap \Gamma\hooklongrightarrow D^{(p)}\times E^{(p)}\yrightarrow{\pr_1} D^{(p)}$
is proper, $\Gamma''^{(p)}$ induces the action on cohomology
$(\Gamma''^{(p)})^*\colon H_c^q(D^{(p)},\Lambda)\longrightarrow H_c^q(E^{(p)},\Lambda)$.
By the projection formula, these two maps are equal.
Now we state the functoriality result:

\subsubsection{Theorem}\label{Thm:functoriality-corr}
Let the notation be the same as above. Then we have a map of spectral sequences
\[
 \xymatrix{%
 E^{p,q}_1=\bigoplus_{i\ge\max(0,-p)}H^{q-2i}_c\bigl(D^{(p+2i)},\Lambda(-i)\bigr)
 \ar@{=>}[r]\ar[d]^-{\oplus (\Gamma''^{(p+2i)})^*}&
 H^{p+q}_c(X_s,R\psi_X\Lambda)\ar[d]^-{\Gamma^*}\\
 E'^{p,q}_1=\bigoplus_{i\ge\max(0,-p)}H^{q-2i}_c\bigl(E^{(p+2i)},\Lambda(-i)\bigr)
 \ar@{=>}[r]&
 H^{p+q}_c(Y_s,R\psi_Y\Lambda)\lefteqn{.}
 }
\]

\begin{prf}
 This follows from Corollary \ref{Cor:special-nearby-corr}, Proposition \ref{Prop:functoriality-pull-back},
 Proposition \ref{Prop:functoriality-cup-product}, and Proposition \ref{Prop:functoriality-push-forward}.
\end{prf}

\section{$\boldsymbol{\ell}$-independence of nearby cycle cohomology}\label{section:nearby-l-ind}
\subsection{$\boldsymbol{\ell}$-independence of nearby cycle cohomology}
\subsubsection{}
Let $K$ be a complete discrete valuation field with finite residue field $F=\F_q$.
We denote the ring of integers of $K$ by $\mathcal{O}_K$ and the characteristic of $F$ by $p$. 
Fix a separable closure $\overline{K}$ of $K$ and let $\overline{F}$ be the residue field of
the integral closure of $\mathcal{O}_K$ in $\overline{K}$, which is an algebraic closure of $F$.
We denote by $G_K$ (resp.\ $G_F$) the Galois group $\Gal(\overline{K}/K)$ (resp.\ $\Gal(\overline{F}/F)$).
We denote by $\Fr_q$ the geometric Frobenius element (the inverse of the $q$th power map)
in $G_F$. The Weil group $W_K$ of $K$ is defined as the inverse image of the subgroup
$\langle \Fr_q\rangle\subset G_F$ by the canonical map $G_K\longrightarrow G_F$.
For $\sigma\in W_K$, let $n(\sigma)$ be the integer such that the image of $\sigma$ in $G_F$
is $\Fr_q^{n(\sigma)}$. Put $W_K^+=\{\sigma\in W_K\mid n(\sigma)\ge 0\}$.

Put $S=\Spec \mathcal{O}_K$. For an $S$-scheme $X$, we denote its special fiber, geometric special fiber,
generic fiber, geometric generic fiber by $X_F$, $X_{\overline{F}}$, $X_K$, $X_{\overline{K}}$ respectively.

Let $\ell$ be a prime number distinct from $p$.

\subsubsection{}
The main result in this section is the following theorem:

\subsubsection{Theorem}\label{Thm:nearby-l-ind}
Let $X$ be a flat arithmetic $S$-scheme with purely $d$-dimensional smooth generic fiber,
and $\Gamma\subset X\times_S X$ a closed subscheme with purely $d$-dimensional generic fiber.
Assume that the composite $\Gamma\hooklongrightarrow X\times_S X\yrightarrow{\pr_1}X$ is proper.
Then for any $\sigma \in W_K^+$, the number
\[
 \Tr\bigl(\Gamma^*\circ \sigma_*;H^*_c(X_{\overline{F}},R\psi\Q_\ell)\bigr)=
 \sum_{i=0}^{2d}(-1)^i\Tr\bigl(\Gamma^*\circ \sigma_*;H^i_c(X_{\overline{F}},R\psi\Q_\ell)\bigr)
\] 
lies in $\Z[1/p]$ and is independent of $\ell$.

\subsubsection{}
First we treat the case where $X$ is strictly semistable over $S$.
We need a slight generalization of the above theorem in this case:

\subsubsection[][\cf {\cite[Lemma 3.2]{TSaito}}]{Lemma}\label{Lemma:nearby-l-ind-ssred}
Let $L$ be a finite quasi-Galois extension of $K$ and put $S'=\Spec \mathcal{O}_L$.
We denote the residue field of $L$ by $E$.
Let $X$ be a strictly semistable $S'$-scheme which is purely of relative dimension $d$.
Take any $\sigma\in W_K^+$. Fix an embedding $\overline{K}\hooklongrightarrow \overline{L}$ and extend $\sigma$ uniquely
to an automorphism of $\overline{L}$.
We put $X^\sigma=X\times_{\mathcal{O}_L\nearrow\sigma}\mathcal{O}_L$.
Let $\Gamma\subset X^\sigma\times_{S'} X$ be a closed subscheme with purely $d$-dimensional
generic fiber satisfying that the composite $\Gamma\longrightarrow X^\sigma\times_{S'} X\yrightarrow{\pr_1}X^{\sigma}$
is proper. Then the number
\[
 \Tr\bigl(\Gamma^*\circ \sigma_*;H^*_c(X_{\overline{E}},R\psi\Q_\ell)\bigr)
\] 
lies in $\Z[1/p]$ and is independent of $\ell$.

\begin{prf}
 We denote the irreducible components of $X_E$ by $D_1,\ldots,D_m$ as usual.
 Then the irreducible components of $X^{\sigma}_E$ are $D^\sigma_1,\ldots,D^\sigma_m$.
 We define $\Gamma''^{(s)}\in \CH_{d-s}((D^{\sigma(s)}\times D^{(s)})\cap \Gamma)$ for each $s$
 as in (\ref{functoriality-corr-notation}).
 Then by Theorem \ref{Thm:functoriality-corr} we have the following map of spectral sequences:
 \[
 \xymatrix{%
 E'^{s,t}_1=\bigoplus_{i\ge\max(0,-s)}H^{t-2i}_c\bigl(D^{\sigma(s+2i)}_{\overline{E}},\Q_\ell(-i)\bigr)
 \ar@{=>}[r]\ar[d]^-{\oplus (\Gamma''^{(s+2i)})^*}&
 H^{s+t}_c(X^\sigma_{\overline{E}},R\psi_{X^\sigma}\Q_\ell)\ar[d]^-{\Gamma^*}\\
 E^{s,t}_1=\bigoplus_{i\ge\max(0,-s)}H^{t-2i}_c\bigl(D^{(s+2i)}_{\overline{E}},\Q_\ell(-i)\bigr)
 \ar@{=>}[r]&
 H^{s+t}_c(X_{\overline{E}},R\psi_X\Q_\ell)\lefteqn{.}
 }
\]
 On the other hand, we have the map of spectral sequences induced by $\sigma$:
 \[
 \xymatrix{%
 E^{s,t}_1=\bigoplus_{i\ge\max(0,-s)}H^{t-2i}_c\bigl(D^{(s+2i)}_{\overline{E}},\Q_\ell(-i)\bigr)
 \ar@{=>}[r]\ar[d]^-{\overline{\sigma}_*}&
 H^{s+t}_c(X_{\overline{E}},R\psi_X\Q_\ell)\ar[d]^-{\sigma_*}\\
 E'^{s,t}_1=\bigoplus_{i\ge\max(0,-s)}H^{t-2i}_c\bigl(D^{\sigma(s+2i)}_{\overline{E}},\Q_\ell(-i)\bigr)
 \ar@{=>}[r]&
 H^{s+t}_c(X^\sigma_{\overline{E}},R\psi_{X^\sigma}\Q_\ell)\lefteqn{,}
 }
\]
 where $\overline{\sigma}$ denotes the image of $\sigma$ in $G_E$.
 Let $\sigma_{\mathrm{geom}}^{(s)}\colon D^{\sigma(s)}_{\overline{E}}\longrightarrow D^{(s)}_{\overline{E}}$
 be the composition $\varphi^{f\cdot n(\sigma)}\circ \overline{\sigma}^*$, where $\varphi$ denotes the absolute Frobenius
 morphism and $f$ is the integer satisfying $q=p^f$.
 This is a proper morphism over $\overline{E}$. Since $\varphi$ induces the identity map on \'etale cohomology,
 we have $\overline{\sigma}_*=\sigma_{\mathrm{geom}}^{(s)*}$. Therefore we obtain the endomorphism of 
 a spectral sequence
  \[
 \xymatrix{%
 E^{s,t}_1=\bigoplus_{i\ge\max(0,-s)}H^{t-2i}_c\bigl(D^{(s+2i)}_{\overline{E}},\Q_\ell(-i)\bigr)
 \ar@{=>}[r]\ar[d]^-{(\Gamma''^{(s+2i)})^*\circ{\sigma_{\mathrm{geom}}^{(s+2i)*}}}&
 H^{s+t}_c(X_{\overline{E}},R\psi_X\Q_\ell)\ar[d]^-{\Gamma^*\circ\sigma_*}\\
 E^{s,t}_1=\bigoplus_{i\ge\max(0,-s)}H^{t-2i}_c\bigl(D^{(s+2i)}_{\overline{E}},\Q_\ell(-i)\bigr)
 \ar@{=>}[r]&
 H^{s+t}_c(X_{\overline{E}},R\psi_X\Q_\ell)\lefteqn{.}
 }
\]
 Denote by $\Gamma'''^{(s)}\in \CH_{d-s}((D_{\overline{E}}^{\sigma(s)}\times D_{\overline{E}}^{(s)})\cap (\sigma^{(s)}_{\mathrm{geom}}\times\id)(\Gamma_{\overline{E}}))$
 the image of $\Gamma''^{(s)}$ under the map
 \begin{align*}
 &\CH_{d-s}\bigl((D^{\sigma(s)}\times D^{(s)})\cap \Gamma\bigr)\longrightarrow 
 \CH_{d-s}\bigl((D_{\overline{E}}^{\sigma(s)}\times D_{\overline{E}}^{(s)})\cap\Gamma_{\overline{E}}\bigr)\\
 &\qquad\qquad\yrightarrow{(\sigma_{\mathrm{geom}}^{(s)}\times\id)_*}
 \CH_{d-s}\Bigl((D_{\overline{E}}^{\sigma(s)}\times D_{\overline{E}}^{(s)})\cap (\sigma^{(s)}_{\mathrm{geom}}\times\id)(\Gamma_{\overline{E}})\Bigr).
 \end{align*}
 Then $(\Gamma'''^{(s)})^*=(\Gamma''^{(s)})^*\circ \sigma_{\mathrm{geom}}^{(s)*}$ holds.
 Thus we have equalities
 \begin{align*}
  \Tr\bigl(\Gamma^*\circ \sigma_*;H^*_c(X_{\overline{E}},R\psi\Q_\ell)\bigr)
  &=\sum_{s}\sum_{i\ge \max(0,-s)}(-1)^s\Tr\Bigl((\Gamma''^{(s+2i)})^*\circ\sigma_{\mathrm{geom}}^{(s+2i)*};H_c^*\bigl(D_{\overline{E}}^{(s+2i)},\Q_\ell(-i)\bigr)\Bigr)\\
  &=\sum_{s}\sum_{i\ge \max(0,-s)}(-1)^sq^{n(\sigma)i}\Tr\bigl((\Gamma'''^{(s+2i)})^*;H_c^*(D_{\overline{E}}^{(s+2i)},\Q_\ell)\bigr)\\
  &=\sum_{s}(-1)^s\cdot\frac{q^{n(\sigma)(s+1)}-1}{q^{n(\sigma)}-1}
  \Tr\bigl((\Gamma'''^{(s)})^*;H_c^*(D_{\overline{E}}^{(s)},\Q_\ell)\bigr).
 \end{align*}
 By Theorem \ref{Thm:finite-l-ind}, the number
 $\Tr\bigl((\Gamma'''^{(s)})^*;H_c^*(D_{\overline{E}}^{(s)},\Q_\ell)\bigr)$ lies in $\Z[1/p]$ and is independent of $\ell$.
 Therefore $\Tr\bigl(\Gamma^*\circ \sigma_*;H^*_c(X_{\overline{E}},R\psi\Q_\ell)\bigr)$ lies in $\Z[1/p]$ and 
 is independent of $\ell$.
\end{prf}

\subsubsection{}\label{subsub:reduce-alteration}
Next we reduce Theorem \ref{Thm:nearby-l-ind} to Lemma \ref{Lemma:nearby-l-ind-ssred} by de Jong's alteration (\cite{deJong}).
We may assume that $X$ is connected. Since $X$ is flat over $S$ with smooth generic fiber, it is 
irreducible and reduced. Therefore by \cite[Theorem 6.5]{deJong} and \cite[Lemma 1.11]{TSaito}, we have a
finite quasi-Galois extension $L$ of $K$, a scheme $Y$ which is strictly semistable over $\mathcal{O}_L$, and
a proper surjective generically finite $S$-morphism $f\colon Y\longrightarrow X$.
Put $S'=\Spec \mathcal{O}_L$ and denote the residue field of $L$ by $E$ as in the proof of
Lemma \ref{Lemma:nearby-l-ind-ssred}.
Let $K'$ be the inseparable closure of $K$ in $L$. Then we have a canonical isomorphism
$H^i_c(X'_{\overline{F}},R\psi_{X'}\Q_\ell)\cong H^i_c(X_{\overline{F}},R\psi_X\Q_\ell)$, where
$X'=X\otimes_{\mathcal{O}_K}\mathcal{O}_{K'}$. Moreover, if we fix an embedding 
$\overline{K}\hooklongrightarrow \overline{K'}$, the isomorphism above is compatible with an isomorphism
$W_{K'}\yrightarrow{\sim} W_K$. Therefore by replacing $K$ and $X$ by $K'$ and $X'$ respectively,
we may assume that the extension $L/K$ is separable. 

We denote by $Y'$ the scheme $Y$ considered as an $S$-scheme. Take any $\sigma\in W_K^+$.
By Lemma \ref{Lemma:nearby-refined-pullback-corr}, we have the commutative diagram below:
\[
\xymatrix{%
H^i_c(Y'_{\overline{F}},R\psi_{Y'}\Q_\ell)\ar[r]^-{\sigma_*}\ar[d]^-{f_*}&
H^i_c(Y'^\sigma_{\overline{F}},R\psi_{Y'^\sigma}\Q_\ell)\ar[r]^-{\Gamma'^*}\ar[d]^-{f^\sigma_*}&
H^i_c(Y'_{\overline{F}},R\psi_{Y'}\Q_\ell)\\
H^i_c(X_{\overline{F}},R\psi_X\Q_\ell)\ar[r]^-{\sigma_*}&
H^i_c(X_{\overline{F}},R\psi_X\Q_\ell)\ar[r]^-{\Gamma^*}&
H^i_c(X_{\overline{F}},R\psi_X\Q_\ell)\ar[u]_-{f^*}\lefteqn{,}
}
\]
where $\Gamma'\in Z_d((f^{\sigma}\times f)^{-1}(\Gamma))$ is an element satisfying
$\Gamma'_K=(f_K^{\sigma}\times f_K)^![\Gamma_K]\in \CH_d((f_K^{\sigma}\times f_K)^{-1}(\Gamma_K))$,
as in Lemma \ref{Lemma:nearby-refined-pullback-corr}.
Together with Lemma \ref{Lemma:nearby-alt}, we have 
\[
 \Tr\bigl(\Gamma'^*\circ \sigma_*; H^i_c(Y'_{\overline{F}},R\psi_{Y'}\Q_\ell)\bigr)
=\deg f\cdot \Tr\bigl(\Gamma^*\circ \sigma_*; H^i_c(X_{\overline{F}},R\psi_X\Q_\ell)\bigr)
\]
as in the proof of \cite[Lemma 3.3]{TSaito}.

Let $h\colon \coprod_{\tau\in\Gal(L/K)}Y^\tau\longrightarrow Y'\otimes_SS'$ be the morphism
induced by the canonical map 
$\mathcal{O}_L\otimes_{\mathcal{O}_K}\mathcal{O}_L\longrightarrow \prod_{\tau\in\Gal(L/K)}\mathcal{O}_L$.
It is finite surjective and induces an isomorphism on generic fibers.
Therefore we have an isomorphism 
$H_c^i(Y'_{\overline{F}},R\psi_{Y'}\Q_\ell)\cong \bigoplus_{\tau\in\Gal(L/K)}H_c^i(Y^{\tau}_{\overline{E}},R\psi_{Y^\tau}\Q_\ell)$.
Under this isomorphism, the map 
$\sigma_*\colon H^i_c(Y'_{\overline{F}},R\psi_{Y'}\Q_\ell)\longrightarrow H^i_c(Y'^{\sigma}_{\overline{F}},R\psi_{Y'^\sigma}\Q_\ell)$ is identified with the direct sum of $\sigma_*\colon H^i_c(Y^{\tau}_{\overline{E}},R\psi_{Y^\tau}\Q_\ell)\longrightarrow H^i_c(Y^{\sigma\tau}_{\overline{E}},R\psi_{Y^{\sigma\tau}}\Q_\ell)$.
For $\tau,\tau'\in\Gal(L/K)$, let 
\[
 \Gamma'_{\tau,\tau'}\in Z_d\bigl((f^\tau\times f^{\tau'})^{-1}(\Gamma)\bigr)
\]
be an element such that $(\Gamma'_{\tau,\tau'})_L=\Gamma'_L\vert_{Y^\tau_L\times Y^{\tau'}_L}$,
where $\Gamma'_L$ is the base change of $\Gamma'_K$ from $K$ to $L$.
By Lemma \ref{Lemma:nearby-refined-pullback-corr} again, the $(\tau,\tau')$-component of the map
\[
 \bigoplus_{\tau\in\Gal(L/K)}H_c^i(Y^{\sigma\tau}_{\overline{E}},R\psi_{Y^{\sigma\tau}}\Q_\ell)\longrightarrow 
\bigoplus_{\tau'\in\Gal(L/K)}H_c^i(Y^{\tau'}_{\overline{E}},R\psi_{Y^{\tau'}}\Q_\ell)
\]
induced by 
$\Gamma'^*\colon H^i_c(Y'^\sigma_{\overline{F}},R\psi_{Y'^\sigma}\Q_\ell)\longrightarrow H^i_c(Y'_{\overline{F}},R\psi_{Y'}\Q_\ell)$
is equal to 
\[
 \Gamma_{\sigma\tau,\tau'}'^*\colon H^i_c(Y^{\sigma\tau}_{\overline{E}},R\psi_{Y^{\sigma\tau}}\Q_\ell)\longrightarrow H^i_c(Y^{\tau'}_{\overline{E}},R\psi_{Y^{\tau'}}\Q_\ell).
\]

Therefore the number
\begin{align*}
\Tr\bigl(\Gamma^*\circ \sigma_*; H^*_c(X_{\overline{F}},R\psi_X\Q_\ell)\bigr)
&=\frac{1}{\deg f}\cdot\Tr\bigl(\Gamma'^*\circ \sigma_*; H^*_c(Y'_{\overline{F}},R\psi_{Y'}\Q_\ell)\bigr)\\
 &=\frac{1}{\deg f}\sum_{\tau\in\Gal(L/K)}\Tr\bigl(\Gamma_{\sigma\tau,\tau}'^*\circ \sigma_*; H^*_c(Y^{\tau}_{\overline{E}},R\psi_{Y^\tau}\Q_\ell)\bigr)
\end{align*}
lies in $\frac{1}{\deg f}\Z[1/p]$ and is independent of $\ell$ by Lemma \ref{Lemma:nearby-l-ind-ssred}.

By the same technique as in \cite[p.~629]{TSaito}, we can derive from the following lemma that the number
$\Tr\bigl(\Gamma^*\circ \sigma_*; H^*_c(X_{\overline{F}},R\psi_X\Q_\ell)\bigr)$ is in $\Z[1/p]$.
Now the proof of Theorem \ref{Thm:nearby-l-ind} is complete.

\subsubsection{Lemma}
Let $K$ be a field of characteristic $0$. Let $a_1,\ldots,a_r$ be distinct elements of $K$ and
$c_1,\ldots,c_r$ non-zero integers. 
Put $s_m=\sum_{i=1}^rc_ia_i^m$ for a non-negative integer $m$.
Assume that there exists an integer $N\ge 1$ such that $Ns_m\in \Z[1/p]$ for every $m\ge 0$.
Then $s_m\in \Z[1/p]$ for every $m\ge 0$.

\begin{prf}
 By \cite[Lemma 2.8]{Kleiman}, $a_i$ is integral over $\Z[1/p]$ for every $i$.
 Therefore every $s_m$ is also integral over $\Z[1/p]$, while it is in $\Q$.
 Since $\Z[1/p]$ is normal, we have $s_m\in \Z[1/p]$.
\end{prf}

\subsubsection{Remark}\label{Remark:nearby-integrality}
{\upshape The result of S.~Bloch and H.~Esnault \cite{Bloch-Esnault} implies that the alternating sum of the trace
 in Theorem \ref{Thm:nearby-l-ind} lies in $\Z$ (\cf Remark \ref{Remark:Bloch-Esnault}).
 For $\Gamma=\Delta_X$ (the diagonal of $X$), the integrality also follows from \cite[Theorem 4.2]{rigid-Weil}.
}

\subsection{$\boldsymbol{\ell}$-independence for stalks of nearby cycles}
\subsubsection{}
In this subsection, we give some results on $\ell$-independence for stalks of nearby cycles.
All of them are immediate consequences of Theorem \ref{Thm:nearby-l-ind}.

\subsubsection{Theorem}\label{Thm:stalk-l-ind}
Let $X$ be a flat arithmetic $S$-scheme with purely $d$-dimensional smooth generic fiber,
and $x\in X_F$ an $F$-rational point. Choose a geometric point $\overline{x}$ lying over $x$.
Then the Weil group $W_K$ acts on the stalk $(R^i\psi_X\Q_\ell)_{\overline{x}}$. For every $\sigma\in W_K^+$, the number
\[
 \Tr\bigl(\sigma_*; (R^*\psi_X\Q_\ell)_{\overline{x}}\bigr)=
\sum_{i=0}^d(-1)^i\Tr\bigl(\sigma_*; (R^i\psi_X\Q_\ell)_{\overline{x}}\bigr)
\]
is an integer which is independent of $\ell$.

\begin{prf}
 Put $U=X\setminus \{x\}$. Then we have the following $W_K$-equivariant exact sequence:
 \[
  \longrightarrow H^i_c(U_{\overline{F}},R\psi_U\Q_\ell)\longrightarrow 
  H^i_c(X_{\overline{F}},R\psi_X\Q_\ell)\longrightarrow (R^i\psi_X\Q_\ell)_{\overline{x}}
 \longrightarrow H^{i+1}_c(U_{\overline{F}},R\psi_U\Q_\ell)\longrightarrow.
 \]
 Therefore we have the equality
 \[
  \Tr\bigl(\sigma_*; (R^*\psi_X\Q_\ell)_{\overline{x}}\bigr)=
  \Tr\bigl(\sigma_*; H^*_c(X_{\overline{F}},R\psi_X\Q_\ell)\bigr)-
  \Tr\bigl(\sigma_*; H^*_c(U_{\overline{F}},R\psi_U\Q_\ell)\bigr).
 \]
 Since each term of the right hand side lies in $\Z[1/p]$ and is independent of $\ell$, 
 so is the left hand side $\Tr\bigl(\sigma_*; (R^*\psi_X\Q_\ell)_{\overline{x}}\bigr)$.

 The integrality follows from Remark \ref{Remark:nearby-integrality} (note that
 since we only use the case $\Gamma=\Delta_X$,
 we do not need the result of S.~Bloch and H.~Esnault).
\end{prf}

\subsubsection{Corollary}\label{Cor:totdim}
Let the notation be the same as in Theorem \ref{Thm:stalk-l-ind}.
Then the integers
\begin{align*}
 \dim_{\Q_\ell}(R^*\psi_X\Q_\ell)_{\overline{x}}
 =\sum_{i=0}^{d}(-1)^i\dim_{\Q_\ell}(R^i\psi_X\Q_\ell)_{\overline{x}},\quad
 \Sw (R^*\psi_X\Q_\ell)_{\overline{x}}=\sum_{i=0}^{d}(-1)^i\Sw(R^i\psi_X\Q_\ell)_{\overline{x}}
\end{align*}
are independent of $\ell$. Here $\Sw$ denotes the Swan conductor.

\begin{prf}
 These are immediate consequences of Theorem \ref{Thm:stalk-l-ind}
 (for the part of the Swan conductor, see \cite[Corollary 2.6]{Ochiai}).
\end{prf}

\subsubsection{Remark}
{\upshape
The above corollary gives a weak evidence of Deligne's conjecture on Milnor numbers
(\cite[Expos\'e XVI, Conjecture 1.9]{SGA7-II}). The statement of the conjecture is the following:
}

\subsubsection{Conjecture}\label{Conj:Deligne}
{\upshape
 Let $K^{\mathrm{ur}}$ be the maximal unramified extension of $K$ and put
 $S^{\mathrm{ur}}=\Spec \mathcal{O}_{K^{\mathrm{ur}}}$.
 Let $X$ be a purely $d$-dimensional flat arithmetic $S^{\mathrm{ur}}$-scheme. Assume that $X$ is regular
 and that the structure morphism $X\longrightarrow S^{\mathrm{ur}}$ is smooth outside a unique closed point
 $x\in X_{\overline{F}}$. Put 
 \begin{align*}
  \dimtot_{\F_\ell}(R^*\phi\,\F_\ell)_x&=\dim_{\F_\ell}(R^*\phi\,\F_\ell)_{x}+\Sw(R^*\phi\,\F_\ell)_{x},\\
  \mu(X/S^{\mathrm{ur}},x)&=\length_{\mathcal{O}_{X,x}}\sExt^1(\Omega_{X/S},\mathcal{O}_X)_x.
 \end{align*}
 Then the equality 
 \[
   \dimtot_{\F_\ell}(R^*\phi\,\F_\ell)_x=\mu(X/S^{\mathrm{ur}},x)	    
 \]
 holds.
 (The original conjecture allows a more general base trait. See \cite{Orgogozo}.)
}

\bigbreak
{\upshape
 This conjecture is solved in the cases below:
 \begin{itemize}
  \item $d=0$ (\cite[Expos\'e XVI, Proposition 1.12]{SGA7-II}).
  \item the point $x$ is an ordinary double point (\cite[Expos\'e XVI, Proposition 1.13]{SGA7-II}).
  \item the characteristic of $K$ is positive (\cite[Expos\'e XVI, Theorem 2.4]{SGA7-II}).
  \item $d=1$ (\cite[Corollaire 0.9]{Orgogozo}).
 \end{itemize}
 Moreover, F.~Orgogozo proved that the conductor formula of Bloch implies the above conjecture
 (\cite[Th\'eor\`eme 0.8]{Orgogozo}).

 Since 
 \begin{align*}
  \dimtot_{\F_\ell}(R^*\phi\,\F_\ell)_{x}
  &=\dim_{\F_\ell}(R^*\psi\F_\ell)_{x}+\Sw(R^*\psi\F_\ell)_{x}-1\\
  &=\dim_{\Q_\ell}(R^*\psi\Q_\ell)_{x}+\Sw(R^*\psi\Q_\ell)_{x}-1
 \end{align*}
 (the last equality follows from the universal coefficient theorem), from Corollary \ref{Cor:totdim} 
 we see that the left hand side
 of the equality in Conjecture \ref{Conj:Deligne} is independent of $\ell$, while the right hand side is obviously
 independent of $\ell$.
 }

\subsection{$\boldsymbol{\ell}$-independence for open schemes over local fields}
\subsubsection{}
In this subsection, we consider an analogue of \cite[Theorem 0.1]{TSaito} for open schemes
over local fields. 

\subsubsection{Definition}
{\upshape Let $X$ be an arithmetic $S$-scheme and $H\subset X$ a closed subscheme of $X$. 
We may write $H=H_h\cup H'$, where $H'$ is contained in the special fiber of $X$ and $H_h\longrightarrow S$ is flat.
The pair $(X,H)$ is called a {\slshape strictly semistable pair} if the following conditions hold
(\cf \cite[6.3]{deJong}):
\begin{itemize}
 \item $X$ is strictly semistable over $S$.
 \item $H$ is a strict normal crossing divisor of $X$.
 \item Let $H_i$ ($i\in I$) be the irreducible components of $H_h$. For each $J\subset I$, the scheme
       $H_J=\bigcap_{i\in J}H_i$ is a union of schemes which are strictly semistable over $S$.
\end{itemize} 
Moreover if $H$ is flat over $S$ (namely, $H=H_h$), we call $(X,H)$ a {\slshape horizontal strictly semistable pair}.
For a strictly semistable pair $(X,H)$, the pair $(X,H_h)$ is a horizontal strictly semistable pair.
}

\subsubsection{Lemma}\label{Lemma:open-vs-nearby}
Let $(X,H)$ be a horizontal strictly semistable pair over $S$.
Put $U=X\setminus H$ and denote the canonical open immersion $U\hooklongrightarrow X$ by $j$.
Then the canonical morphism 
\[
 j_{\overline{F}!}R\psi_U\Q_\ell\longrightarrow R\psi_X(j_{\overline{K}!}\Q_\ell)
\]
is an isomorphism. In particular, if $X$ is proper over $S$, we have an isomorphism
$H^i_c(U_{\overline{F}},R\psi_U\Q_\ell)\cong H^i_c(U_{\overline{K}},\Q_\ell)$.

\begin{prf}
 Since the problem is \'etale local, we may assume that 
 \[
  X=\Spec\mathcal{O}_K[T_1,\ldots,T_n]/(T_{r+1}\cdots T_s-\pi),\quad H=V(T_1\cdots T_r)\subset X,
 \]
 where $\pi$ is a uniformizer of $K$ (\cf \cite[1.5 (d)]{Illusie-Dwork}).
 Put $X_1=\Spec\mathcal{O}_K[T_{r+1},\ldots,T_n]/(T_{r+1}\cdots T_s-\pi)$.
 Then $(X,H)\cong (\A^r_S\times_SX_1,Z\times_SX_1)$, where
 $Z\subset \A^r_S$ is the divisor defined by $T_1\cdots T_r=0$.
 By the K\"unneth formula for $R\psi$ (\cite[Th\'eor\`eme 4.7]{Illusie-autour}), 
 we may reduce the lemma to the case $(X,H)=(\A^r_S,Z)$. 
 This case is treated in \cite[Expos\'e XIII, Proposition 2.1.9]{SGA7-II}.
\end{prf}

\subsubsection{}
The following proposition is an analogue of Lemma \ref{Lemma:nearby-l-ind-ssred}:

\subsubsection{Proposition}\label{Prop:open-l-ind-sspair}
Let $L$ be a finite quasi-Galois extension of $K$ and put $S'=\Spec \mathcal{O}_L$.
We denote the residue field of $L$ by $E$.
Let $X$ be an arithmetic $S'$-scheme with purely $d$-dimensional generic fiber
and assume that there exists a compactification $X\hooklongrightarrow \overline{X}$ over $S'$ such that
$(\overline{X},\overline{X}\setminus X)$ is a strictly semistable pair over $S'$.
Take any $\sigma\in W_K^+$. Fix an embedding $\overline{K}\hooklongrightarrow \overline{L}$ and extend $\sigma$ uniquely
to an automorphism of $\overline{L}$.
We put $X^\sigma=X\times_{\mathcal{O}_L\nearrow\sigma}\mathcal{O}_L$.
Let $\Gamma\subset X^\sigma\times_{S'} X$ be a closed subscheme with purely $d$-dimensional
generic fiber satisfying that the composite $\Gamma\longrightarrow X^\sigma\times_{S'} X\yrightarrow{\pr_1}X^{\sigma}$
is proper. Then the number
\[
 \Tr\bigl(\Gamma^*_{L}\circ \sigma_*;H^*_c(X_{\overline{L}},\Q_\ell)\bigr)
\] 
lies in $\Z[1/p]$ and is independent of $\ell$.

\begin{prf}
 By Lemma \ref{Lemma:open-vs-nearby}, 
 $H^i_c(X_{\overline{L}},\Q_\ell)\cong H^i_c(X_{\overline{E}},R\psi_X\Q_\ell)$ and
 $H^i_c(X^\sigma_{\overline{L}},\Q_\ell)\cong H^i_c(X^\sigma_{\overline{E}},R\psi_{X^\sigma}\Q_\ell)$ hold.
 Furthermore we can easily see that the map 
 $\Gamma_L^*\colon H^i_c(X^\sigma_{\overline{L}},\Q_\ell)\longrightarrow H^i_c(X_{\overline{L}},\Q_\ell)$ 
 corresponds to the map 
 $\Gamma^*\colon H^i_c(X^\sigma_{\overline{E}},R\psi_{X^\sigma}\Q_\ell)\longrightarrow H^i_c(X_{\overline{E}},R\psi_X\Q_\ell)$
 (\cf Proposition \ref{Prop:partial-nearby-corr}). Thus the number
 \[
 \Tr\bigl(\Gamma_L^*\circ \sigma_*;H^*_c(X_{\overline{L}},\Q_\ell)\bigr)=
 \Tr\bigl(\Gamma^*\circ \sigma_*;H^*_c(X_{\overline{E}},R\psi_X\Q_\ell)\bigr)
 \]
 lies in $\Z[1/p]$ and is independent of $\ell$ by Lemma \ref{Lemma:nearby-l-ind-ssred}.
\end{prf}

\subsubsection{}\label{subsub:horizontal-ss}
Let $X$ be a scheme which is smooth and separated of finite type over $K$. 
Take a compactification $X\hooklongrightarrow Z$ over $S$.
Namely, $Z$ is a scheme which is proper and flat over $S$, containing $X$ as an open subscheme.
Put $Y=Z\setminus X$. By de Jong's alteration (\cite[Theorem 6.5]{deJong}), there exist a finite extension $L$ of $K$, 
an arithmetic $\mathcal{O}_L$-scheme $W$, a proper surjective generically finite $S$-morphism
$f\colon W\longrightarrow Z$ such that $(W,f^{-1}(Y))$ is a strictly semistable pair over $S'=\Spec \mathcal{O}_L$.
Let $H$ be the horizontal part $f^{-1}(Y)_h$ of $f^{-1}(Y)$. Then $(W,H)$ is a horizontal strictly
semistable pair over $S'$
such that $(W\setminus H)_K\longrightarrow X$ is a proper surjective generically finite $K$-morphism.

By the lemma below, we can take $L$ as a quasi-Galois extension of $K$.

\subsubsection{Lemma}
Let $(X,H)$ be a horizontal strictly semistable pair over $S$. Let $L$ be a finite extension of $K$ and put 
$S'=\Spec\mathcal{O}_L$. Then there exists a blow-up $\pi\colon X'\longrightarrow X\times_SS'$ whose center is contained
in the special fiber such that $(X',\pi^{-1}(H))$ is a horizontal strictly semistable pair over $S'$.

\begin{prf}
 We may take the same blow-up as in \cite[Lemma 1.11]{TSaito}.
\end{prf}

\subsubsection{Theorem}\label{Thm:open-l-ind}
Let $X$ be a purely $d$-dimensional scheme which is smooth and separated of finite type over $K$,
and $\Gamma\subset X\times X$ a purely $d$-dimensional closed subscheme such that the composite
$\Gamma\hooklongrightarrow X\times X\yrightarrow{\pr_1} X$ is proper.
Let $Z$, $L$, $(W,H)$, $f\colon W\longrightarrow Z$ be as in (\ref{subsub:horizontal-ss}) (we take $L$
as a quasi-Galois extension of $K$). Put $U=W\setminus H$ and write $g\colon U_L\longrightarrow X$ for the restriction
of $f$. 
Assume that the composite 
$\overline{(g\times g)^{-1}(\Gamma)}\hooklongrightarrow U\times_{S'}U\yrightarrow{\pr_1}U$ is proper
($\overline{(g\times g)^{-1}(\Gamma)}$ denotes the closure of $(g\times g)^{-1}(\Gamma)\subset U_L\times U_L$
in $U\times_{S'}U$).
Then for any $\sigma \in W_K^+$, the number
\[
 \Tr\bigl(\Gamma^*\circ \sigma_*;H^*_c(X_{\overline{K}},\Q_\ell)\bigr)=
 \sum_{i=0}^{2d}(-1)^i\Tr\bigl(\Gamma^*\circ \sigma_*;H^i_c(X_{\overline{K}},\Q_\ell)\bigr)
\] 
lies in $\Z[1/p]$ and is independent of $\ell$.

\begin{prf}
 As in (\ref{subsub:reduce-alteration}), we may assume that the extension $L/K$ is separable.
 Put $V=U_L$.
 We denote by $V'$ the scheme $V$ considered as a scheme over $K$. We have 
 $V'_L\cong \coprod_{\tau\in\Gal(L/K)}V^{\tau}$. Take any $\sigma\in W_K^+$ and
 put $\Gamma'=(g^{\sigma}\times g)^![\Gamma]\in \CH_d((g^{\sigma}\times g)^{-1}(\Gamma))$. 
 For $\tau,\tau'\in\Gal(L/K)$, put $\Gamma'_{\tau,\tau'}=\Gamma_L'\vert_{V^\tau\times V^{\tau'}}$,
 where $\Gamma'_L$ is the base change of $\Gamma'$ from $K$ to $L$.
 As in (\ref{subsub:reduce-alteration}), we have the equality
 \[
  \Tr\bigl(\Gamma^*\circ \sigma_*;H^*_c(X_{\overline{K}},\Q_\ell)\bigr)
 =\frac{1}{\deg f}\sum_{\tau\in\Gal(L/K)}\Tr\bigl(\Gamma'^*_{\sigma\tau,\tau}\circ \sigma_*;H^*_c(V_{\overline{L}},\Q_\ell)\bigr).
 \]
 By the assumption, for each $\tau,\tau'\in \Gal(L/K)$, there exists a cycle 
 $\overline{\Gamma}'_{\tau,\tau'}\in Z_d(U^\tau\times_{S'}U^{\tau'})$
 such that $(\overline{\Gamma}'_{\tau,\tau'})_L=\Gamma'_{\tau,\tau'}$ and the composite 
 $\lvert \overline{\Gamma}'_{\tau,\tau'}\rvert\hooklongrightarrow U^\tau\times_{S'}U^{\tau'}\yrightarrow{\pr_1}U^{\tau}$
 is proper. Therefore we may reduce our theorem to Proposition \ref{Prop:open-l-ind-sspair}.
\end{prf}

\section{$\boldsymbol{\ell}$-independence for rigid spaces}\label{section:rigid-l-ind}
Let the notation be the same as in the previous section.
We consider rigid spaces over a complete discrete valuation field $K$
as adic spaces locally of finite type over $\Spa (K,\mathcal{O}_K)$ (\cf \cite{Huber-generalization}).
We denote a scheme by an ordinary italic letter such as $X$,
a formal scheme by a calligraphic letter such as $\mathcal{X}$, and 
a rigid space by a sans serif letter such as $\X$.
For a scheme $X$ over $S=\Spec \mathcal{O}_K$, we denote the completion of $X$ along its special fiber
by $X^{\wedge}$. For a formal scheme $\mathcal{X}$ over $\Spf \mathcal{O}_K$,
we write $\mathcal{X}^\rig$ for its Raynaud generic fiber.
It is the analytic adic space $d(\mathcal{X})$ in \cite[1.9]{Huber-book}.

\subsection{Smooth case}
\subsubsection{}
In this subsection, we prove our main theorem for smooth rigid spaces.
We derive the following consequence from the result in the previous section.

\subsubsection{Corollary}\label{Cor:algebraizable-l-ind}
Let $X$ be an arithmetic $S$-scheme with smooth generic fiber
and $\X$ the rigid space $(X^{\wedge})^\rig$.
Then for every $\sigma\in W_K^+$, the number
\[
 \Tr\bigl(\sigma_*;H^*_c(\X_{\overline{K}},\Q_\ell)\bigr)
 =\sum_{i=0}^{2\dim \X}(-1)^i\Tr\bigl(\sigma_*;H^i_c(\X_{\overline{K}},\Q_\ell)\bigr)
\]
is an integer which is independent of $\ell$.

\begin{prf}
 We may assume that $X$ is connected and flat over $S$.
 We have a $W_K$-equivariant isomorphism $H^i_c(\X_{\overline{K}},\Q_\ell)\cong H^i_c(X_{\overline{F}},R\psi_X\Q_\ell)$
 (\cite[Theorem 5.7.6]{Huber-book}).
 Applying Theorem \ref{Thm:nearby-l-ind} to $\Gamma=X\stackrel{\Delta_X}{\hooklongrightarrow}X\times_SX$, 
 we see that for every $\sigma\in W_K^+$ the number
 \[
 \Tr\bigl(\sigma_*;H^*_c(\X_{\overline{K}},\Q_\ell)\bigr)=\Tr\bigl(\sigma_*;H^*_c(X_{\overline{F}},R\psi_X\Q_\ell)\bigr)
 \]
 lies in $\Z[1/p]$ and is independent of $\ell$. On the other hand, we know that every eigenvalue of 
 the action of $\sigma\in W_K^+$ on $H^i_c(\X_{\overline{K}},\Q_\ell)$ is an algebraic integer
 (\cite[Theorem 4.2]{rigid-Weil}).
 Therefore the rational number $\Tr\bigl(\sigma_*;H^*_c(\X_{\overline{K}},\Q_\ell)\bigr)$
 is an algebraic integer, i.e., an integer.
\end{prf}

\subsubsection{Definition}
{\upshape A formal scheme $\mathcal{X}$ of finite type over $\mathcal{O}_K$ is said to be of {\slshape type (SA)}
(smoothly algebraizable) if there exists an arithmetic $S$-scheme $X$ with smooth generic fiber
such that $\mathcal{X}\cong X^{\wedge}$.
A rigid space $\X$ over $K$ is said to be of {\slshape type (SA)} if there exists a formal scheme $\mathcal{X}$ 
of type (SA) over $\mathcal{O}_K$ such that $\X\cong \mathcal{X}^\rig$.
}

\subsubsection{Lemma}\label{Lem:adm-blowup-algebraizable}
Let $X$ be an arithmetic $S$-scheme with smooth generic fiber. Then the following hold:
\begin{enumerate}
 \item Every admissible blow-up of $X^\wedge$ is of type (SA).
 \item Every open formal subscheme of $X^\wedge$ is of type (SA).
\end{enumerate}

\begin{prf}
 \begin{enumerate}
  \item Take a uniformizer $\pi$ of $K$.
	Let $\mathcal{I}$ be an open ideal of $\mathcal{O}_{X^\wedge}$. Since the topology of $\mathcal{O}_{X^\wedge}$
	is the $\pi$-adic topology and $X^\wedge$ is noetherian, there exists an integer $n$ satisfying
	$\pi^n\mathcal{O}_{X^\wedge}\subset \mathcal{I}$. Denote by $\mathcal{I}'$ the unique ideal of $\mathcal{O}_X$
	containing $\pi^n\mathcal{O}_X$ such that $\mathcal{I}/\pi^n\mathcal{O}_{X^\wedge}
	=\mathcal{I'}/\pi^n\mathcal{O}_X$.
	It is clear that $\mathcal{I'}\mathcal{O}_{X^\wedge}$ coincides with $\mathcal{I}$.
	Then the admissible blow-up of $X^\wedge$ by $\mathcal{I}$ is equal to the $\pi$-adic completion of
	the scheme $X'$ obtained by the blow-up of $X$ by $\mathcal{I}'$.
	The generic fiber of $X'$ is obviously smooth.
  \item We can identify the underlying topological space of $X^\wedge$ with that of $X_F$.
	Let $\mathcal{U}$ be an open formal subscheme of $X^\wedge$.
	Then, $U=X\setminus (X_F\setminus \mathcal{U})$ is an arithmetic open subscheme of $X$
	satisfying $U_F=\mathcal{U}$ as topological spaces. Then the generic fiber of $U$ is smooth and 
	$U^\wedge=\mathcal{U}$.
 \end{enumerate}
\end{prf}

\subsubsection{Corollary}\label{Cor:algebraizable-open}
Let $\X=(X^{\wedge})^\rig$ be a rigid space of type (SA) over $K$. 
Then every quasi-compact open subspace $\mathsf{U}$ of\, $\X$ is of type (SA).

\begin{prf}
 Since $\mathsf{U}$ is quasi-compact, there exist an admissible blow-up 
 $\mathcal{Y}\longrightarrow X^{\wedge}$ and
 an open formal subscheme $\mathcal{U}\subset \mathcal{Y}$ such that $\mathsf{U}=\mathcal{U}^\rig$
 (\cite[Lemma 4.4]{BL1}).
 By Lemma \ref{Lem:adm-blowup-algebraizable}, $\mathcal{Y}$ and $\mathcal{U}$ are of type (SA).
 This completes the proof.
\end{prf}

\subsubsection{Theorem}\label{Thm:smooth-l-ind}
Let $\X$ be a quasi-compact separated rigid space which is smooth over $K$. 
Then for every $\sigma\in W_K^+$, the number
\[
 \Tr\bigl(\sigma_*;H^*_c(\X_{\overline{K}},\Q_\ell)\bigr)
\]
is an integer which is independent of $\ell$. 

\begin{prf}
 By \cite[Corollary 2.5]{rigid-Weil}, there exists a finite open covering $\{\mathsf{U}_i\}_{1\le i\le m}$ of\, $\X$
 consisting of rigid spaces of type (SA). Corollary \ref{Cor:algebraizable-open} ensures that each intersection
 $\mathsf{U}_{i_1}\cap \cdots\cap \mathsf{U}_{i_n}$ is of type (SA). Thus by Corollary \ref{Cor:algebraizable-l-ind},
 for every $\sigma\in W_K^+$, the number
 \[
 \Tr\Bigl(\sigma_*;H^*_c\bigl((\mathsf{U}_{i_1}\cap\cdots\cap \mathsf{U}_{i_n})_{\overline{K}},\Q_\ell\bigr)\Bigr)
 \]
 is an integer which is independent of $\ell$. On the other hand, we have the spectral sequence below:
 \[
 E_1^{-s,t}=\bigoplus_{1\le i_1<\cdots<i_s\le m}H^t_c\bigl((\mathsf{U}_{i_1}\cap \cdots\cap \mathsf{U}_{i_s})_{\overline{K}},\Q_\ell\bigr)
 \Longrightarrow H^{-s+t}_c(\X_{\overline{K}},\Q_\ell).
 \]
 Therefore the number
\[
 \Tr\bigl(\sigma_*;H^*_c(\X_{\overline{K}},\Q_\ell)\bigr)
 =\sum_{s=1}^{m}(-1)^s\kern-10pt\sum_{1\le i_1<\cdots<i_s\le m}
 \Tr\Bigl(\sigma_*;H^*_c\bigl((\mathsf{U}_{i_1}\cap\cdots\cap \mathsf{U}_{i_s})_{\overline{K}},\Q_\ell\bigr)\Bigr)
\] 
 is also an integer which is independent of $\ell$. 
\end{prf}

\subsubsection{}
From now on we consider ordinary cohomology. First we establish the analogous result as in 
\cite[Theorem 4.2]{rigid-Weil}.

\subsubsection{Theorem}\label{Thm:smooth-Weil-ord}
Let $\X$ be a quasi-compact separated rigid space which is smooth over $K$. 
Then for every $\sigma\in W_K^+$, every eigenvalue $\alpha\in \overline{\Q}_\ell$ of its action on
$H^i_c(\X_{\overline{K}},\Q_\ell)$
is an algebraic integer. Moreover, there exists a non-negative integer $m$ such that for any isomorphism
$\iota\colon \overline{\Q}_\ell\yrightarrow{\sim}\C$, the absolute value $\lvert \iota(\alpha)\rvert$ is
equal to $q^{n(\sigma)\cdot m/2}$.

\begin{prf}
 We may assume that $\X$ is connected. Put $d=\dim \X$. By the Poincar\'e duality (\cite[Corollary 7.5.6]{Huber-book}),
 every eigenvalue $\alpha$ of $\sigma_*$ on $H^i(\X_{\overline{K}},\Q_\ell)$ is of the form
 $q^{n(\sigma)\cdot d}/\beta$, where $\beta$ is an eigenvalue of $\sigma_*$ on
 $H^{2d-i}_c(\X_{\overline{K}},\Q_\ell)$. Therefore $\alpha$ is an algebraic number and
 there exists an integer $m$ such that for any isomorphism $\iota\colon \overline{\Q}_\ell\yrightarrow{\sim}\C$,
 the absolute value $\lvert \iota(\alpha)\rvert$ is equal to $q^{n(\sigma)\cdot m/2}$.

 Thus we have only to show that $\alpha$ is an algebraic integer. By the same method as in \cite[\S 4]{rigid-Weil},
 we can reduce the theorem to the case $\X=(X^{\wedge})^{\rig}$, where $X$ is strictly semistable scheme over $S$.
 Furthermore by using an analogue of weight spectral sequence, we may reduce the claim to the lemma below
 (\cf \cite[proof of Proposition 4.7]{rigid-Weil}).
\end{prf}

\subsubsection{Lemma}
Let $X$ be a scheme separated of finite type over $\F_q$. Then every eigenvalue of the action of $\Fr_q$ on 
$H^i(X_{\overline{\F}_q},\Q_\ell)$ is an algebraic integer
(here $\Fr_q\in \Gal(\overline{\F}_q/\F_q)$ is the geometric Frobenius element).

\begin{prf}
 We may assume that $X$ is irreducible. By de Jong's alteration \cite{deJong}, we may assume that there exist
 a proper smooth scheme $\overline{X}$ and a strict normal crossing divisor $D$ of $\overline{X}$ such that 
 $X=\overline{X}\setminus D$. 
 Let $D_1,\ldots,D_m$ be the irreducible components of $D$. Put $D_I=\bigcap_{i\in I}D_i$ for 
 $I\subset \{1,\ldots,m\}$ ($D_I=\overline{X}$ for $I=\varnothing$)
 and $D^{(k)}=\coprod_{I\subset \{1,\ldots,m\},\#I=k}D_I$.
 By the spectral sequence 
 \[
  E_1^{-k,n+k}=H^{n-k}\bigl(D^{(k)}_{\overline{\F}_q},\Q_\ell(-k)\bigr)\Longrightarrow H^n(X_{\overline{\F}_q},\Q_\ell),
 \]
 the eigenvalue $\alpha$ occurs as an eigenvalue of $\Fr_{D^{(k)}}^*$ on 
 $H^{n-k}(D^{(k)}_{\overline{\F}_q},\Q_\ell(-k))$ for some $n$, $k$.
 Since $D^{(k)}$ is proper smooth over $\F_q$,
 \cite[Corollaire 3.3.3]{Weil2} ensures that $\alpha$ is integral over $\Z$.
\end{prf}

\subsubsection{Theorem}\label{Thm:smooth-l-ind-ord}
Let $\X$ be a quasi-compact separated rigid space which is smooth over $K$. 
Then for every $\sigma\in W_K^+$, the number
\[
 \Tr\bigl(\sigma_*;H^*(\X_{\overline{K}},\Q_\ell)\bigr)
 =\sum_{i=0}^{2\dim\X}(-1)^i\Tr\bigl(\sigma_*;H^i(\X_{\overline{K}},\Q_\ell)\bigr)
\]
is an integer which is independent of $\ell$. 

\begin{prf}
 By Theorem \ref{Thm:smooth-Weil-ord}, it is sufficient to show that the number
 $\Tr\bigl(\sigma_*;H^*(\X_{\overline{K}},\Q_\ell)\bigr)$ is a rational number which is independent of $\ell$.
 We may assume that $\X$ is connected. Put $d=\dim \X$.
 Let $\alpha_{\ell,i,1},\ldots,\alpha_{\ell,i,m_i}$ be the eigenvalues of $\sigma_*$ on 
 $H^i_c(\X_{\overline{K}},\Q_\ell)$. Then by the Poincar\'e duality, 
 the eigenvalues of $\sigma_*$ on $H^{2d-i}(\X_{\overline{K}},\Q_\ell)$ are
 $q^{n(\sigma)\cdot d}\alpha_{\ell,i,1}^{-1},\ldots,q^{n(\sigma)\cdot d}\alpha_{\ell,i,m_i}^{-1}$. Therefore
 it is sufficient to prove that the number $\sum_{i=0}^{2\dim\X}\sum_{j=1}^{m_i}(-1)^i\alpha_{\ell,i,j}^{-1}$
 is a rational number which is independent of $\ell$. 
 For every non-negative integer $k$, by applying Theorem \ref{Thm:smooth-l-ind} to $\sigma^k\in W_K^+$, we 
 can see that the number
 $\sum_{i=0}^{2\dim\X}\sum_{j=1}^{m_i}(-1)^i\alpha_{\ell,i,j}^k$ is a rational number which is independent of $\ell$.
 As in the proof of Lemma \ref{Lemma:vandermond}, we may conclude that 
 the number $\sum_{i=0}^{2\dim\X}\sum_{j=1}^{m_i}(-1)^i\alpha_{\ell,i,j}^{-1}$
 is a rational number which is independent of $\ell$. 
\end{prf}

\subsection{General case}
\subsubsection{}
In this subsection, we prove our main theorem for general rigid spaces over local fields of characteristic $0$.
We need the following continuity theorem of R.~Huber, which is stronger than
\cite[Proposition 2.1 (iv)]{Huber-comparison} (\cf \cite[Theorem 5.3]{rigid-Weil}).

\subsubsection{Theorem}\label{Thm:strong-continuity}
Assume that the characteristic of $K$ is equal to $0$.
Let $\X$ be a quasi-compact separated rigid space over $K$ and $\mathsf{Z}$ a closed analytic subspace of\, $\X$.
Write $\mathsf{U}$ for $\X\setminus \mathsf{Z}$.
Then for every pair of prime numbers $\ell$, $\ell'$ which do not divide $q$,
there exists a quasi-compact open subspace $\mathsf{U}'$ of $\mathsf{U}$ such that 
the canonical maps $H^i_c(\mathsf{U}'_{\overline{K}},\Z_\ell)\longrightarrow H^i_c(\mathsf{U}_{\overline{K}},\Z_\ell)$
and $H^i_c(\mathsf{U}'_{\overline{K}},\Z_{\ell'})\longrightarrow H^i_c(\mathsf{U}_{\overline{K}},\Z_{\ell'})$
are isomorphisms for every $i$.

\begin{prf}
 This is due to \cite[(II) in the proof of Theorem 3.3]{Huber-comparison}. We briefly recall the argument there.
 By \cite[Corollary 2.7]{Huber-JAG-!}, there exists $\varepsilon_0>0$ such that for every $0<\varepsilon<\varepsilon_0$
 the canonical map $H^i_c(\mathsf{U}(\varepsilon)_{\overline{K}},\Z/\ell\Z)
 \longrightarrow H^i_c(\mathsf{U}_{\overline{K}},\Z/\ell\Z)$ is an isomorphism.
 Here we write $\mathsf{U(\varepsilon)}$ for $P(\varepsilon)$ in \cite[2.6]{Huber-JAG-!}.
 By the long exact sequence of cohomology groups derived from the short exact sequence of sheaves
 \[
  0\longrightarrow \Z/\ell\Z\yrightarrow{\times \ell^n} \Z/\ell^{n+1}\Z\longrightarrow \Z/\ell^n\Z
 \longrightarrow 0,
 \]
 we see inductively that the canonical map 
 $H^i_c(\mathsf{U}(\varepsilon)_{\overline{K}},\Z/\ell^n\Z)\longrightarrow H^i_c(\mathsf{U}_{\overline{K}},\Z/\ell^n\Z)$
 is an isomorphism for every $0<\varepsilon<\varepsilon_0$ and $n$. 
 In the same way, there exists $\varepsilon_1>0$ such that for every $0<\varepsilon<\varepsilon_1$ and $n$ the canonical map
 $H^i_c(\mathsf{U}(\varepsilon)_{\overline{K}},\Z/\ell'^n\Z)\longrightarrow H^i_c(\mathsf{U}_{\overline{K}},\Z/\ell'^n\Z)$
 is an isomorphism. Put $\varepsilon_2=\min\{\varepsilon_0,\varepsilon_1\}$ and $\mathsf{U}'=\mathsf{U}(\varepsilon_2)$.
 Then $\mathsf{U}'$ is quasi-compact and both of the canonical maps
 \[
  H^i_c(\mathsf{U}'_{\overline{K}},\Z/\ell^n\Z)\longrightarrow H^i_c(\mathsf{U}_{\overline{K}},\Z/\ell^n\Z),
 \quad 
  H^i_c(\mathsf{U}'_{\overline{K}},\Z/\ell'^n\Z)\longrightarrow H^i_c(\mathsf{U}_{\overline{K}},\Z/\ell'^n\Z)
 \]
 are isomorphisms.

 On the other hand we have the canonical isomorphisms 
 \begin{align*}
 &\varprojlim_{n}H^i_c(\mathsf{U}'_{\overline{K}},\Z/\ell^n\Z)\cong H^i_c(\mathsf{U}'_{\overline{K}},\Z_\ell),&
 &\varprojlim_{n}H^i_c(\mathsf{U}_{\overline{K}},\Z/\ell^n\Z)\cong H^i_c(\mathsf{U}_{\overline{K}},\Z_\ell),\\
 &\varprojlim_{n}H^i_c(\mathsf{U}'_{\overline{K}},\Z/\ell'^n\Z)\cong H^i_c(\mathsf{U}'_{\overline{K}},\Z_{\ell'}),&
 &\varprojlim_{n}H^i_c(\mathsf{U}_{\overline{K}},\Z/\ell'^n\Z)\cong H^i_c(\mathsf{U}_{\overline{K}},\Z_{\ell'})  
 \end{align*}
 (\cite[Theorem 3.1 and Theorem 3.3]{Huber-comparison}).
 Therefore the canonical homomorphisms
 \[
  H^i_c(\mathsf{U}'_{\overline{K}},\Z_\ell)\longrightarrow H^i_c(\mathsf{U}_{\overline{K}},\Z_\ell),\quad
  H^i_c(\mathsf{U}'_{\overline{K}},\Z_{\ell'})\longrightarrow H^i_c(\mathsf{U}_{\overline{K}},\Z_{\ell'})
 \]
 are isomorphisms.
\end{prf}

\subsubsection{Theorem}\label{Thm:general-l-ind}
Assume that the characteristic of $K$ is equal to $0$.
Let $\X$ be a quasi-compact separated rigid space over $K$. Then for every $\sigma\in W_K^+$,
the number
\[
 \Tr\bigl(\sigma_*;H^*_c(\X_{\overline{K}},\Q_\ell)\bigr)
 =\sum_{i=0}^{2\dim\X}(-1)^i\Tr\bigl(\sigma_*;H^i_c(\X_{\overline{K}},\Q_\ell)\bigr)
\]
is an integer which is independent of $\ell$.

\begin{prf}
 Let $\ell$ and $\ell'$ be prime numbers which do not divide $q$ and $\sigma\in W^+_K$.
 We prove by induction on $\dim \X$ that the numbers
 \[
 \Tr\bigl(\sigma_*;H^*_c(\X_{\overline{K}},\Q_\ell)\bigr),\qquad
 \Tr\bigl(\sigma_*;H^*_c(\X_{\overline{K}},\Q_{\ell'})\bigr) 
 \]
 are integers and are equal.
 We may assume that $\X$ is reduced.
 Let $\mathsf{Z}$ be the singular locus of $\X$.
 It is a closed analytic subspace whose dimension is strictly less than $\dim \X$.
 Thus we have only to show our claim on $H^i_c(\mathsf{U}_{\overline{K}},\Q_\ell)$ and
 $H^i_c(\mathsf{U}_{\overline{K}},\Q_{\ell'})$,
 where $\mathsf{U}=\mathsf{X}\setminus \mathsf{Z}$.
 Take an open subspace $\mathsf{U}'\subset \mathsf{U}$ as in Theorem \ref{Thm:strong-continuity}.
 Then we have the isomorphisms
 \[
 H^i_c(\mathsf{U}'_{\overline{K}},\Q_\ell)\yrightarrow {\sim} H^i_c(\mathsf{U}_{\overline{K}},\Q_\ell),\qquad
 H^i_c(\mathsf{U}'_{\overline{K}},\Q_{\ell'})\yrightarrow {\sim} H^i_c(\mathsf{U}_{\overline{K}},\Q_{\ell'})
 \]
 by Theorem \ref{Thm:strong-continuity}.
 Therefore by Theorem \ref{Thm:smooth-l-ind} the numbers
 \[
 \Tr\bigl(\sigma_*;H^*_c(\mathsf{U}_{\overline{K}},\Q_\ell)\bigr),\qquad
 \Tr\bigl(\sigma_*;H^*_c(\mathsf{U}_{\overline{K}},\Q_{\ell'})\bigr) 
 \]
 are integers and are equal. This completes the proof.
\end{prf}

\end{document}